%% file: JCAM_revision.tex
\documentclass[preprint,11pt]{elsarticle}

\usepackage{graphicx}
\usepackage[utf8]{inputenc}
\usepackage{multirow}
\usepackage{comment}
\usepackage{mdwlist}
\usepackage{amsfonts}
\usepackage{ulem}
\usepackage{amsmath} 
\usepackage{amsthm}
\usepackage{amssymb}
\usepackage{titlesec}
\usepackage[utf8]{inputenc}

\usepackage{lipsum}
\usepackage{amsfonts}
\usepackage{graphicx}
\usepackage{epstopdf}

\usepackage{fancyhdr}

\usepackage{algorithm}
\usepackage{algorithmic}

\usepackage{hyperref}      
\usepackage{url} 
\usepackage{titlesec}
\usepackage{indentfirst}
\usepackage{booktabs}
\usepackage{verbatim}
\usepackage{color}
\usepackage{latexsym}
\usepackage{comment}
\usepackage{amscd}
\usepackage{esvect}
\usepackage{graphicx}
\usepackage{upgreek}
\usepackage{caption}   
\usepackage[a4paper]{geometry}
\usepackage{mathrsfs}
\usepackage[numbers]{natbib}
\usepackage{float}
\usepackage{lipsum}
\usepackage[export]{adjustbox} 
\usepackage{subfiles}
\usepackage{subfigure}
\usepackage{bm}
\usepackage{epstopdf}

\newcommand{\bx}{\boldsymbol{x}}

\newcommand{\bv}{\boldsymbol{v}}

\newcommand{\blue}{\textcolor{black}}

\begin{document}

\begin{frontmatter}
\author{Zhenyi Zhu}
\ead{zyzhu@math.cuhk.edu.hk}
\author{Yuchen Huang}
\ead{ychuang@math.cuhk.edu.hk}
\author{Liu Liu}
\ead{lliu@math.cuhk.edu.hk}
\address{Department of Mathematics, The Chinese University of Hong Kong}
\title{PhysicsSolver: Transformer-Enhanced Physics-Informed Neural Networks for Forward and Forecasting Problems in Partial Differential Equations}
\begin{abstract}
Time-dependent partial differential equations are a significant class of equations that describe the evolution of various physical phenomena over time. One of the open problems in scientific computing is predicting the behaviour of the solution outside the given temporal region. Most traditional numerical methods are applied to a given time-space region and can only accurately approximate the solution of the given region. To address this problem, many deep learning-based methods, basically data-driven and data-free approaches, have been developed to solve these problems. However, most data-driven methods require a large amount of data, which consumes significant computational resources and fails to utilize all the necessary information embedded underlying the partial differential equations (PDEs). Moreover, data-free approaches such as Physics-Informed Neural Networks (PINNs) may not be that ideal in practice, as traditional PINNs, which primarily rely on multilayer perceptrons (MLPs) and convolutional neural networks (CNNs), tend to overlook the crucial temporal dependencies inherent in real-world physical systems. We propose a method denoted by \textbf{PhysicsSolver} that merges the strengths of two approaches: data-free methods that can learn the intrinsic properties of physical systems without using data, and data-driven methods, which are effective at making predictions. 
Extensive numerical experiments have demonstrated the efficiency and robustness of our proposed method. We provide the code at \href{https://github.com/PhysicsSolver/PhysicsSolver}{https://github.com/PhysicsSolver}.
\end{abstract}
\end{frontmatter}

\section{Introduction}

Solving partial differential equations is crucial for many real-world applications, including weather forecasting, industrial design, and material analysis. Numerical methods for solving PDEs, such as the finite element method \cite{bathe2006finite} and the pseudospectral method \cite{fornberg1998practical}, have been extensively studied in science and engineering, they often face high computational costs, especially when dealing with complex or high-dimensional PDEs. Mostly, these traditional numerical methods are mainly designed to obtain the solution in a given region. Thus, they can not forecast the solution outside the given temporal region with high accuracy. To predict the solution outside the given region, they often rely on extrapolation methods \cite{sidi2003practical} such as linear, polynomial, and conic extrapolation, among others.

Recently, data-free machine learning techniques, particularly PINNs, have shown promise in transforming scientific computing by offering faster solutions that can complement traditional methods \cite{raissi2019physics}. However, conventional neural networks typically work within fixed, finite-dimensional spaces, which limits their ability to generalize solutions across different discretizations. This limitation highlights the need to develop mesh-invariant neural networks. Also, these methods are often data-free, which can make the prediction less accurate than data-driven methods in modelling some complex systems.

In addition, data-driven deep learning models such as those based on operator learning have emerged as effective tools to solve PDEs. These models leverage their strong capability for non-linear mapping to learn the relationships between inputs and outputs in PDE-related tasks by utilizing the provided training data. As a result, they can provide solutions much more quickly than traditional numerical methods during the inference phase. One of the most popular data-driven models is the Fourier Neural Operator (FNO) method \cite{li2020fourier}, which can efficiently predict higher-resolution solutions based on lower-resolution inputs. However, a significant challenge remains in the requirement for large datasets to train these models effectively. Moreover, unlike PINNs, which integrate such information directly, these methods do not fully leverage the physics information provided by PDEs. As a result, they may fail to capture some of the intrinsic details behind the datasets.

To address the shortcomings of the above three approaches, we propose a model based on the mixture of transformer and PINNs denoted as \textbf{PhysicsSolver}. Specifically, we first perform data engineering, where we introduce the Interpolation Pseudo Sequential Grids Generator module to convert pointwise spatiotemporal inputs into temporal sequences. Then, we apply Monte Carlo sampling using the Halton sequence, which allows us to sample data points that are relatively sparse yet more representative, providing better datasets for later training. To learn intrinsic physical patterns and their evolution over time, we propose the PhysicsSolver module based on transformers incorporating the Physics-Attention module. This module can encode both the input dataset and the grid points, and then generate the output through a decoder. During the training process, we develop a robust algorithm utilizing both the equation information and sparse training data to optimize the model simultaneously. Finally, we conduct extensive experiments on several one-dimensional and high-dimensional PDEs, where PhysicsSolver achieves consistent state-of-the-art with sharp relative gains. 

The main contributions of this work are summarized as the following: 
\vspace{-5pt}
\begin{itemize}
    \item First, we formulate general forecasting problems in PDEs, then study the capability of different models to predict the future physics states, which is generally a challenging problem in numerical analysis.
    \item Beyond prior methods, we propose a model named PhysicsSolver that blends the advantages of PINNs and the transformer, enabling our model to better learn the evolution of physics states. We found that our proposed PhysicsSolver can capture the trend of intrinsic physical phenomena more accurately and be able to better forecast solutions outside the given temporal region, than extrapolation and other deep learning-based methods.
    \item We develop some data engineering modules such as Interpolation Pseudo Sequential Grid Generator (IPSGG) and Halton Sequence Generation (HSG) to tailor the learning process of PhysicsSolver.
    \item PhysicsSolver achieves robust and superior performance in all benchmark problems for both the forward and forecasting problems in PDEs.
\end{itemize}

The rest of the paper is organized as follows. In Section \ref{sec: problem_formualtion}, we first introduce
the problem formulation of forward problems, inverse problems and forecasting problems in time-dependent PDEs. In Section \ref{sec: related_works}, we review the related works on various methods for solving PDEs. In Section \ref{sec: methodology}, 
we provide details of the methodology and learning schemes for our proposed \textbf{PhysicsSolver} method. In Section \ref{sec: numerical}, a series of numerical experiments are shown to validate the effectiveness of our method, where both forward and forecasting problems are investigated. Finally, the conclusion and future work will be discussed in the last section.

\section{Problem Formulation} 
\label{sec: problem_formualtion}
Let $\mathcal{T}$ be the given temporal space, $\mathcal{D}$ be a bounded open set in $\mathbb{R}^N$ (the spatial space), with boundary $\gamma$ a $(N-1)$-dimensional manifold (the boundary space), and $\Omega^{\bv}$ the velocity space. Considering the following PDEs with spatial input $\bx$, velocity input $\bv$, and temporal input $t$ fitting the following abstraction:
\begin{equation}\label{PDEsys2}
    \begin{cases}
    \mathcal{L} \big[ u(t,\boldsymbol{x},\boldsymbol{v},\boldsymbol{\beta})\big] = f(t,\boldsymbol{x},\boldsymbol{v}),\:\: \forall (t, \boldsymbol{x},\boldsymbol{v}) \in \mathcal{T}\times \mathcal{D}\times\Omega^{\boldsymbol{v}},\\
     \mathcal{B} \big[ u(t,\boldsymbol{x},\boldsymbol{v},\boldsymbol{\beta})\big] = g(t,\boldsymbol{x},\boldsymbol{v}),\:\: \forall (t, \boldsymbol{x}, \boldsymbol{v}) \in \mathcal{T}\times \gamma\times\Omega^{\boldsymbol{v}},\\ \mathcal{I}\big[u(0,\boldsymbol{x},\boldsymbol{v},\boldsymbol{\beta})\big] = h(\boldsymbol{x},\boldsymbol{v}),\:\: \forall (\boldsymbol{x},\boldsymbol{v}) \in \mathcal{D}\times\Omega^{\boldsymbol{v}}.\\
    \end{cases}
\end{equation}

where $u$ is the PDE's solution,   $\boldsymbol{\beta} \in \mathcal{H}^{d}$ is the physical parameter of the PDEs and $\mathcal{H}$ is the parameter space, $d$ is the dimension of physical parameters in the PDEs system. $\mathcal{L}$ is a (possibly nonlinear) differential operator that regularizes the behaviour of the system, $\mathcal{B}$ describes the boundary conditions and $\mathcal{I}$ describes the initial conditions in general. Specifically, $(t,\boldsymbol{x},\boldsymbol{v})\in\mathcal{T}\times\mathcal{D}\times \Omega^{\bv}$ are residual points, $(t, \boldsymbol{x}, \boldsymbol{v}) \in \mathcal{T}\times \gamma\times\Omega^{\boldsymbol{v}}$ are boundary points, and $ (0,\boldsymbol{x},\boldsymbol{v}) \in \mathcal{T}\times\gamma\times\Omega^{\boldsymbol{v}}$ are initial points. $f,g,h$ are source terms, (possibly hybrid) boundary values and initial values, respectively.

There have been two main problems in numerical analysis for PDEs: The forward problem and the inverse problem. The forecasting problem is also an important problem but has not been widely investigated or formulated. Next we briefly introduce them.

\textbf{Forward Problem.} In the forward problem, we would like to approximate the true solution $u$ on $n$ given grids $\{t_i,\boldsymbol{x}_i,\boldsymbol{v}_i\}_{i=0}^{n}$ using traditional numerical schemes such as the finite difference method,  or using the approximated solution $u^{\text{NN}}$ represented by the deep neural networks.

\textbf{Inverse Problem.} In the inverse problem \cite{aster2018parameter}, we need to infer the unknown physical parameter $\boldsymbol{\beta}$ of the PDEs. Different from the forward problem, we need to additionally use synthetic dataset $\{u(t_j,\boldsymbol{x}_j,\boldsymbol{v}_j\})_{j=0}^{n_{\text{data}}}$, where $n_{\text{data}}$ is the volume of datasets available for the given PDEs.

\textbf{Forecasting Problem.} In the forecasting problem, we would like to use the learned model to forecast the solution in an unknown time zone. Denote our learned model be $\mathcal{F}(t,\boldsymbol{x},\boldsymbol{v},\theta)$, then the predicted solution is $\mathcal{F} (\tilde{t},\boldsymbol{x},\boldsymbol{v},\theta)$, where $\tilde{t}$ is the unknown time stamps different from the times stamps $t$ in the grids. 

\section{Related Works}
\label{sec: related_works}
\textbf{Traditional Numerical Methods.} As a fundamental scientific problem, obtaining analytical solutions for partial differential equations is often difficult. Consequently, PDEs are usually discretized into meshes and solved using numerical methods in practice, such as finite difference methods \cite{sod1978survey}, finite element methods \cite{huebner2001finite}, spectral methods \cite{shen2011spectral}, etc. However, these numerical methods usually take a few hours or even days for complex structures \cite{umetani2018learning} and perform badly on the forecasting problem.

\textbf{Data-free deep learning Method.} One of the most popular data-free deep learning methods is PINNs. This approach formalizes the constraints of PDEs, including the equations themselves, as well as initial and boundary conditions, as the objective functions within deep learning models \cite{raissi2019physics,yu2018deep}. During the training process, the outputs of these models gradually align with the PDE constraints, allowing them to effectively approximate the solutions to the PDEs. However, this method mainly relies
on convolutional neural networks, ignoring the crucial temporal dependencies inherent in practical physics systems, making it challenging to forecast solutions outside the given training grids. Also, directly using PINNs will fail to learn the solution of PDEs in some complex scenarios. For instance, \cite{jin2023asymptotic,liu2024asymptotic} propose the Asymptotic-Preserving PINNs to solve the difficulties caused by multiscale equations. \blue{Except the methods mentioned above, \cite{fu2024physics} proposed the PIKFNNs, which can improve the performance of the traditional PINNs by embedding PDE information into the activation function of neural networks. }

\textbf{Operator Learning-Based Methods.} 
Operator learning-based methods constitute a significant class within data-driven deep learning approaches, which involve training neural operators to approximate the input-output relationships in tasks governed by partial differential equations. This method can be applied to many physical scenarios, such as predicting future fluid behaviour based on past observations or estimating the internal stress in solid materials \cite{lu2021learning}. Some of the most well-known models in this area include the Fourier Neural Operator \cite{li2020fourier} and its variants, such as those described in  \cite{li2024geometry} and \cite{rahman2022u}. There have been some works \cite{yin2022continuous} focusing on forecasting the dynamics in PDEs. However, they often require a large volume of data and mostly they often ignore utilizing the intrinsic physical mechanism inside the PDEs.

\textbf{Transformer-Based Models.}
The Transformer model \cite{vaswani2017attention} has garnered considerable attention for its capability to capture long-term dependencies, resulting in substantial advancements in natural language processing tasks \cite{kalyan2021ammus}. Additionally, Transformers have been adapted for use in various other fields, including computer vision, speech recognition, and time-series analysis \cite{dong2018speech, han2022survey,wen2022transformers}. So far, there are few researches on the application of transformers in solving PDEs. Recent researchers have used transformer to learn the solution of given PDEs \cite{cao2021choose, wu2024transolver,zhao2023pinnsformer}. However, the combination of PINNs and transformer hasn't been fused that well and forecasting tasks in PDEs have not been studied enough.

\blue{Aside from the topics of solving PDEs, there are several other intriguing subjects worth exploring, such as the training challenges associated with PINNs. For example, \cite{penwarden2023unified} introduces a unified scalable framework designed to address the training difficulties encountered in PINNs and XPINNs, particularly for time-dependent PDEs. Additionally, \cite{guo2023pre} presents a pre-training PINN framework aimed at enhancing the convergence and accuracy of the standard PINNs method, while \cite{shin2020convergence} offers notable theoretical advancements in applying deep learning to solve PDEs.}

\section{Methodology}
\label{sec: methodology}
To address both the forward and forecasting tasks, especially for the forecasting task, we design the innovative model denoted as PhysicsSolver. In this section, we will mainly introduce two parts of the proposed method: model design and learning scheme. In the part of model design, we will present several modules including Interpolation Pseudo Sequential Grids Generator (IPSGG), Halton Sequence Generation (HSG) and Physics-attention (P-Attention). In the part of the learning scheme, we will provide full details such as the loss function and training algorithms.

\subsection{Model Structure}
\begin{figure}[h]
\centerline{\includegraphics[width=1\linewidth]{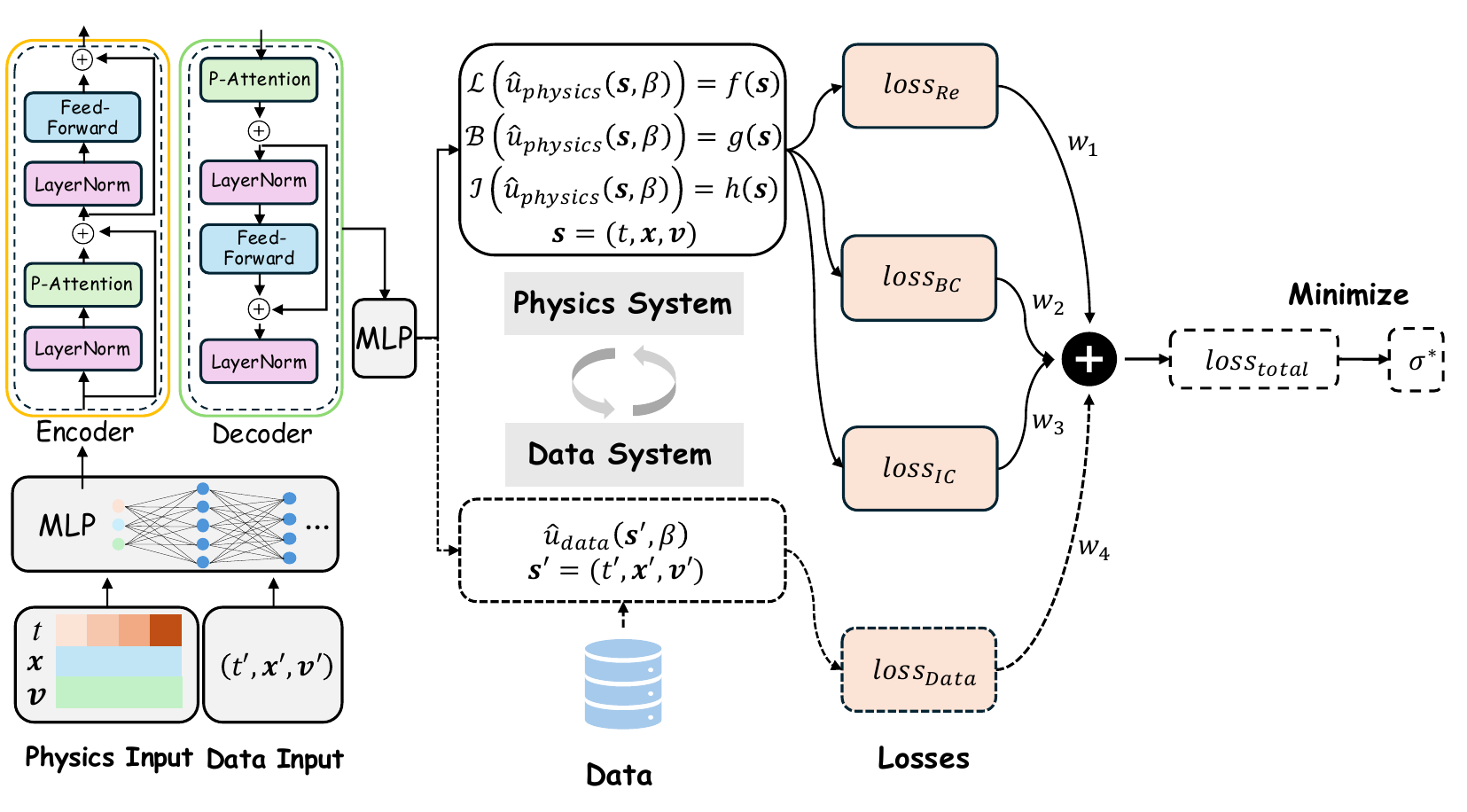}}
    \caption{Architecture of proposed PhysicsSolver.}
    \label{fig:framework}
\end{figure}

We summarize the framework of our proposed PhysicsSolver for the forward and forecasting problem in Figure \ref{fig:framework}. After processing through the IPSGG and HSG modules, the generated grid points are inputted into the Physics-Attention module. During the encoding-decoding process, we can obtain the approximated solutions. One estimated solution $\hat{u}_{\text{physics}}$ will be put into the Physics system for training, while the other estimated solution $\hat{u}_{\text{data}}$ based on Halton sequences will be put into the Data system for later training. We do not utilize any true data from the given PDEs except the solution on sparse grid points generated by HSG module in the training process.
\subsection{Data Generation Engineering}
In this section, we introduce two data generation modules named Interpolation Pseudo Sequential Grids Generator and Halton Sequence Generation designed to prepare the training sequences for PhysicsSolver.
\subsubsection{Interpolation Pseudo Sequential Grids Generator} 
Transformer-based models are designed to capture dependencies within sequential data, whereas traditional PINNs rely on non-sequential data fixed on specific grids for neural network input. Thus, to integrate PINNs with Transformer-based models, we must convert the pointwise spatiotemporal inputs into temporal sequences. At the same time, the generated sequential points should be independent of the given grid points. 
So we design the Interpolation Pseudo Sequential Grids Generator, which can generate sequential grids for the given $p$-th $(t_p,\boldsymbol{x}_p,\boldsymbol{v}_p)$ in the following way:
\begin{equation} \label{IPSGG}
(t_p,\boldsymbol{x}_p,\boldsymbol{v}_p)\xrightarrow{\text{IPSGG}} \{(t_p,\boldsymbol{x}_p,\boldsymbol{v}_p), (t_p+\gamma\Delta t,\boldsymbol{x}_p,\boldsymbol{v}_p), \ldots, (t_p+(k-1)\gamma\Delta t,\boldsymbol{x}_p,\boldsymbol{v}_p)\},
\end{equation}
where $k \in \mathbb{Z}, (k-1)\gamma \not\in \mathbb{Z}$, with $\mathbb{Z}$ representing the set of integers. The parameters $k$ and $\Delta t$ are hyperparameters that determine the number of steps the given grid should move forward and the duration of each step, while $\gamma$ is the control parameter that ensures that the generated sequential points remain independent of the grids used for the forward and forecasting problems. 
\subsubsection{Halton Sequence Generator}
To incorporate extra information from the true data into our PhysicsSolver, we designed the Halton Sequence Generator, which can help sample relatively sparse but more representative data points from the given region. And the use of Halton sequences improves model accuracy compared with data-free models. The Halton sequences method is a type of quasi-Monte Carlo methods, which is tailored to solve some problems in complex systems such as high-dimensional systems as it has multi-dimensional low-discrepancy properties \cite{caflisch1998monte,wang2000randomized}. 
We use the Halton sequences based on the following reasons. Compared with the non-sampling methods, sampling methods can significantly reduce the computational resources. Therefore, the sample size $n_{\text{data}}$ is usually much smaller than the grid size $n_{\text{physics}}$ used in the physics system. Also, quasi-Monte Carlo methods can enhance the accuracy of Monte Carlo approximations based on random sampling \cite{jank2005quasi,shapiro2003monte}.
The generation of the Halton sequences is as follows: Suppose the grid space for a forward problem and a forecasting problem defined above is $\mathcal{X} \in \mathbb{R}^d$, where $d$ is the dimension of the grid points. We use coprime numbers as the bases to construct the Halton sequences from grid space $\mathcal{X'} \in \mathbb{R}^d$. In practice, $\mathcal{X'}$ is typically the complement of $\mathcal{X}$.

\subsection{Physics Informed Neural Networks}
The standard Physics-informed Neural Networks (PINNs) \cite{raissi2019physics} integrate the residual of the underlying PDE, as well as initial and boundary conditions or other relevant physical properties, into the loss function. These networks utilize automatic differentiation to embed physical laws alongside data-driven terms within the loss function, evaluated at specific points in a computational domain.

We begin by reviewing the standard PINNs approach within our problem context. For forward problems, our objective is to approximate the solution of the specified PDEs \eqref{PDEsys2}. Consider a neural network with $L$ layers, the input layer takes $(t,\boldsymbol{x},\boldsymbol{v})$, and the output from the final layer is represented as $u^{\text{NN}}(t,\boldsymbol{x},\boldsymbol{v};m,w,b)$ or, more succinctly, $u^{\text{NN}}_{\theta}$, where $\theta$ denotes the neural network parameters. The relation between the $l$-th and $(l+1)$-th layers ($l=1,2, \cdots, L-1$) is given by
\begin{equation}\label{two-layer}
    n_j^{(l+1)}=\sum_{i=1}^{m_{l}} w_{j i}^{(l+1)} \sigma_{l} (n_{i}^{l})+b_{j}^{(l+1)},
\end{equation} 
where $m=\left(m_{0}, m_{1}, m_{2}, \dots, m_{L-1}\right)$, $w=\left\{w_{j i}^{(k)}\right\}_{i, j, k=1}^{m_{k-1}, m_{k}, L}$ and $b=\left\{b_{j}^{(k)}\right\}_{j=1, k=1}^{m_{k}, L}$. More specifically, 
\begin{itemize}
    \item $n_i^l$: the $i$-th neuron in the $l$-th layer 
    \item $\sigma_l$: the activation function in the $l$-th layer
    \item $w_{j i}^{(l+1)}$: the weight between the $i$-th neuron in the $l$-th layer and the $j$-th neuron in the $(l+1)$-th layer
    \item $b_j^{(l+1)}$: the bias of the $j$-th neuron in the $(l+1)$-th layer
    \item $m_l$: the number of neurons in the $l$-th layer. 
\end{itemize}
Let $\hat{u}$ be neural network approximations, in the traditional PINNs method, we use the following loss function to obtain the constraints on $\hat{u}$.
\begin{equation}
\begin{aligned}
\mathcal{L}_\text{PINNs} :=& w_r \int_{\mathcal{T}} \int_{\mathcal{D}} \int_{\Omega^{\bv}} \Big| \mathcal{L}[\hat{u}(t, \boldsymbol{x}, \boldsymbol{v})]-f(t, \boldsymbol{x}, \boldsymbol{v})\Big|^2  \mathrm{d}\bv \mathrm{d}\bx \mathrm{d}t\\
&+ w_b \int_{\mathcal{T}} \int_\gamma\int_{\Omega^{\boldsymbol{v}}} \Big| \mathcal{B}[\hat{u}(t, \boldsymbol{x}, \boldsymbol{v})]-g(t, \boldsymbol{x}, \boldsymbol{v})\Big|^2  \mathrm{d}\boldsymbol{v}\mathrm{d}\boldsymbol{x} \mathrm{d}t\\
&+ w_i  \int_{\mathcal{D}} \int_{\Omega^{\bv}} \Big| \mathcal{I}[\hat{u}(0, \boldsymbol{x}, \boldsymbol{v})]-h(0, \boldsymbol{x}, \boldsymbol{v})\Big|^2 \mathrm{d}\bv \mathrm{d}\bx. 
\end{aligned}
\label{eq:pinn_loss}
\end{equation}
Where $\mathcal{T}$ is the temporal space, $\mathcal{D}$ is the spatial space, $\mathcal{\gamma}$ is the boundary space, and $\Omega^{\bv}$ is the velocity space. $w_{r}, w_{b}, w_{i}$ are residual weight, boundary weight and initial value weight, respectively. The empirical loss $\mathcal{L}_\text{PINNs\_empirical}$ for equation \eqref{eq:pinn_loss} can be found in equation \eqref{empirical_loss_appendix} from Appendix. In practice, the neural network is trained by minimizing an empirical loss function.

\subsection{Physics-Attention}
The basic structure of PhysicsSolver is the encoder-decoder architecture, which is similar to the Transformer \cite{vaswani2017attention}. In PhysicsSolver, the encoder consists of one self-attention layer called Physics-Attention (P-Attention) and a feedforward layer, where the P-Attention mainly serves to process the mixed embedding containing both physics input information and data input information. The decoder maintains a similar structure to the encoder. Also, both the encoder and the decoder utilize the layer normalization scheme. The encoder is used to process the mixed embedding of both physics input and data input processed by MLP. Specifically, the physics and data inputs are first passed through one MLP, generating two distinct embeddings. These embeddings are then combined into a mixed representation, which is further processed by the encoder.

Intuitively, the self-attention mechanism of P-Attention enables the model to learn the physical dependencies across all spatiotemporal information. At the same time, it can effectively learn the mapping from the input data representation to the output solution, which significantly aids in the prediction task. This capability allows the model to capture more information than traditional PINNs, which mainly focus on approximating the solution of the current state. It also surpasses certain transformer-based models like PINNsformer, which also do not utilize data as effectively.

Suppose we have obtained the mixed embedding $\mathcal{Y} \in \mathcal{R}^{l \times d}$ from the inputs. Here $l$ is the number of the selected grids and $d=2N+1$ is the input dimension for each mesh grid. $\mathcal{Y}$ is first transformed into $\tilde{\mathcal{Y}} \in \mathcal{R}^{l \times d'}$ by the MLP layer, where $d' \geq d$. Then it is transformed into $\mathcal{Y}^{\text{En}}\in \mathcal{R}^{l \times d'}$ by the encoder, where 
\begin{equation} \label{encoder}
    \mathcal{Y}^{\text{En}} = \text{Ln}(\text{Feed}(\text{Attn}(\text{Ln}(\tilde{\mathcal{Y}})))),
\end{equation}
Ln is the layer normalization, Feed is the feed forward layer and Attn is the attention operator. Especially, when we perform the attention operation in both stages of encoding and decoding, we define the trainable projection matrices similar to \cite{cao2021choose}. The latent representations $Q,K,V$ are defined as follows:
\begin{align}
    Q:= \tilde{\mathcal{Y}}W^{Q},\, K:= \tilde{\mathcal{Y}}W^{K},\, V:= \tilde{\mathcal{Y}}W^{V},
\end{align}
where $W^{Q}, W^{K}, W^{V} \in \mathcal{R}^{d' \times d'}$ are trainable projection matrices. The attention operation is defined as
\begin{align}
    \mathrm{Attn}(\tilde{\mathcal{Y}}) := \mathrm{Softmax}\bigg(\dfrac{QK^T}{\Vert QK^T\Vert_{l^2}}\bigg)V.
\end{align}
During the decoding process, the encoded embedding $\mathcal{Y}^{\text{En}}$ is transferred into the decoder, and we obtain the decoded embedding as follows:
\begin{equation} \label{decoder}
    \mathcal{Y}^{\text{De}} = \text{Ln}(\text{Feed}(\text{Ln}(\text{Attn}(\mathcal{Y}^{\text{En}})))),
\end{equation}
where $\mathcal{Y}^{\text{De}} \in \mathcal{R}^{l \times d'}$. After processing through the output MLP layer, $\mathcal{Y}^{\text{De}}$ is converted to $\mathcal{Y}^\text{final}$. $\mathcal{Y}^\text{final}$ can be considered as the final representation of the input time sequence, and serves as the approximated solution $\hat{u}$ for the PDE systems. This representation can subsequently be used as the approximated solutions $\hat{u}_{\text{physics}}$ and $\hat{u}_{\text{data}}$ for the physics system and data system respectively.

\subsection{Learning Scheme}
In this section, we present the learning scheme for PhysicsSolver by adapting the traditional PINNs loss \eqref{eq:pinn_loss}. For the physics system, we consider the residual loss $\mathcal{L}_\text{physics\_res}$, the boundary condition loss $\mathcal{L}_\text{physics\_bc}$ and the initial condition loss $\mathcal{L}_\text{physics\_ic}$ as follows:
\begin{equation}
\begin{aligned}
    \mathcal{L}_\text{physics\_res} = &\frac{1}{kN_\text{r}}\sum_{p=1}^{N_\text{r}}\sum_{j=0}^{k-1} \|\mathcal{L}[\hat{u}_{\text{physics}}(t_p+j \gamma \Delta t, \boldsymbol{x}_p,\boldsymbol{v}_p)] - f(t_p+j \gamma \Delta t,\boldsymbol{x}_p, \boldsymbol{v}_p))\|^2,\\
    \mathcal{L}_\text{physics\_bc} = &\frac{1}{kN_\text{b}}\sum_{p=1}^{N_\text{b}}\sum_{j=0}^{k-1} \|\mathcal{B}[\hat{u}_{\text{physics}}(t_p+j \gamma\Delta t, \boldsymbol{x}_p, \boldsymbol{v}_p)] - g(t_p+j \gamma\Delta t, \boldsymbol{x}_p, \boldsymbol{v}_p))\|^2,\\
    \mathcal{L}_\text{physics\_ic} =& \frac{1}{N_\text{i}} \sum_{p=1}^{N_\text{i}} \|\mathcal{I}[\hat{u}_{\text{physics}}(0, \boldsymbol{x}_p,\boldsymbol{v}_p)] - h(0, \boldsymbol{x}_p,\boldsymbol{v}_p)\|^2. \\
\end{aligned} 
\label{eqn:obj}
\end{equation}
Where $N_\text{r}$, $N_\text{b}$, and $N_\text{i}$ are the number of the residual points, boundary condition points and initial condition points, respectively. The total loss for the physics system is as follows:
\begin{equation}\label{loss_physics}
    \mathcal{L}_{\text{physics}} = w_r \mathcal{L}_\text{physics\_res}  + w_b \mathcal{L}_\text{physics\_bc} + w_i \mathcal{L}_\text{physics\_ic}.
\end{equation}

We also introduce the solution for mismatching loss of the data system:

\begin{equation} \label{loss_data}
    \mathcal{L}_{\text{data}} = \sum_{d=1}^{n_{\text{data}}} (\hat{u}_{\text{data}}(t_d,\bx_d, \bv_d) - u(t_d, \bx_d,\bv_d))^2,
\end{equation}
where the grid points $\{(t_d, \bx_d,\bv_d)\}_{d=1}^{n_{\text{data}}}$ are the data points generated by HSG module and $u(t_d, \bx_d, \bv_d)$ is the groundtruth solution. Finally, the full loss of the PhysicsSolver for the forward and forecasting problem is defined as follows:
\begin{equation} \label{loss_total}
    \mathcal{L}_{\text{PhysicsSolver}} = \lambda_1 \mathcal{L}_\text{physics} 
 + \lambda_2 \mathcal{L}_{\text{data}}.
\end{equation}
Here $\lambda_1$ and $\lambda_2$ are the penalty parameters.

\subsection {Algorithm}
We present the algorithm for training the \textbf{PhysicsSolver} in Algorithm \ref{algorithm_training}. Define $C$ to be the feature dimension of mesh grids we use except the temporal dimension. We have \begin{equation}
    C=X\times V,\, \,C_b= X_b\times V_b, \,\,C_i=X_i\times V_i.
\end{equation} Here $T, T_{b}$ are the numbers of temporal points for the interior domain and the boundary; $X, X_{b}, X_{i}$ are the numbers of spatial points for the interior domain, the boundary and initial domain; $V, V_{b}, V_i$ are the numbers of velocity points used for the interior domain, the boundary and initial domain, respectively. 
\begin{algorithm} [h]
\caption{Training procedure for PhysicsSolver}
\begin{algorithmic}[1] \label{algorithm_training}

\STATE \textbf{Initialization:} Initialize the hyperparameter: learning rate $\eta$, parameter $k$, control parameter $\gamma$ in IPSGG. Initialize the network weight parameters PhysicsSolver module $\theta_{\text{PhysicsSolver}}$, which is the synthetical weight for MLP module and encoder-decoder module
\STATE \textbf{Data Generation:} Generate $n_1$ pseudo sequential grids $\mathbf{X_{\text{physics}}}=\{(t_p,\boldsymbol{x}_p,\boldsymbol{v}_p)\}_{\boldsymbol{p}=1}^{n_1}$ by IPSGG module, see equation \eqref{IPSGG}. We have $\mathbf{X_{\text{physics}}}= \mathbf{X}_{1:N_r} \cup \mathbf{X}_{1:N_b} \cup \mathbf{X}_{1:N_i}$. Here  $\mathbf{X}_{1:N_r} \in \mathbb{R}^{T \times C}$, $\mathbf{X}_{1:N_b} \in \mathbb{R}^{T_b \times C_b}$ and $\mathbf{X}_{1:N_i} \in \mathbb{R}^{C_i}$ are interior, boundary and initial grid points respectively. Generate $n_2$ grid points $\mathbf{X}_{\text{data}} = \{(t_d,\boldsymbol{x}_d,\boldsymbol{v}_d)\}_{\boldsymbol{d}=1}^{n_2}$ by HSG module
\STATE \textbf{Input:} Interior grid points $\mathbf{X}_{1:N_r} \in \mathbb{R}^{T \times C}$, boundary grid points $\mathbf{X}_{1:N_b} \in \mathbb{R}^{T_b \times C_b}$, initial grid points $\mathbf{X}_{1:N_i} \in \mathbb{R}^{C_i}$, and data grid points $\mathbf{X}_{\text{data}}$
\FOR{epoch = 1 to \text{max\_epoch}}
\STATE Feed $\mathbf{X}_{1:N_r}$,$\mathbf{X}_{1:N_b}$, $\mathbf{X}_{1:N_i}$, and  $\mathbf{X}_{\text{data}}$ into MLP layer, followed by the encoder in PhysicsSolver. Then obtain the representation $\mathbf{E}_{1:N_r}$, $\mathbf{E}_{1:N_b}$, $\mathbf{E}_{1:N_i}$, and $\mathbf{E}_{\text{data}}$, see equation \eqref{encoder}
\STATE Decode $\mathbf{E}_{1:N_r}$, $\mathbf{E}_{1:N_b}$,$\mathbf{E}_{1:N_i}$, and $\mathbf{E}_{\text{data}}$ to obtain the representation $\mathbf{D}_{1:N_r}$, $\mathbf{D}_{1:N_b}$,$\mathbf{D}_{1:N_i}$, and $\mathbf{D}_{\text{data}}$, see equation \eqref{decoder}. Processed by the MLP layer in PhysicsSolver, we can obtain the approximated solutions $\hat{u}_{\text{physics}}$ and $\hat{u}_{\text{data}}$ for physics system and data system
\STATE Obtain the training loss for physics system $\mathcal{L}_{\text{physics}}$, see equation \eqref{loss_physics}
\STATE Obtain the training loss for data system $\mathcal{L}_{\text{data}}$, see equation \eqref{loss_data}
\STATE Obtain the total loss for PhysicsSolver $\mathcal{L}_{\text{PhysicsSolver}} = \lambda_1\mathcal{L}_\text{physics} 
 + \lambda_2\mathcal{L}_{\text{data}}$, see equation  \eqref{loss_total}
\STATE \textbf{Gradient Update:} $\theta \gets \theta - \eta \nabla_{\theta} \mathcal{L}_{\text{PhysicsSolver}}(\theta)$, where $\theta =  \theta_{\text{PhysicsSolver}}$
\ENDFOR
\end{algorithmic}
\end{algorithm}

\section{Numerical Examples}
\label{sec: numerical}

To show our model's superiority over traditional numerical methods and other deep learning models, we test the performance of different models on four PDE systems, some of the experimental cases are from \cite{raissi2019physics} and \cite{zhao2023pinnsformer}. For each case, we conduct two kinds of numerical experiments of target PDEs: forward problem and prediction problem. For the forward problem, we test PhysicsSolver with deep learning methods such as PINNs and PINNsformer. For the forecasting problem, we test PhysicsSolver with deep learning methods such as PINNs and PINNsformer and traditional numerical methods such as Extrapolation. In the forecasting problem, we consider the single-step prediction, which means we only predict one unknown time step using the trained models.

\textbf{Networks Architecture.} For \textbf{PhysicsSolver} method, we utilize the Transformer backbone with 2 heads, where the hidden size is set to 512 and the embedding size is set to 32. We set the control parameter $\gamma = 1.1$, and the hyperparameter $k=5$. In different experimental setups, the step size $\Delta t$ varies among $1 \times 10^{-2}$, $1 \times 10^{-3}$, and $1 \times 10^{-4}$, and the number of data points generated by HSG $n_{\texttt{data}}$ is set to either 2, 4, 5, or 10. For \textbf{PINNsFormer} method, we use similar settings in \cite{zhao2023pinnsformer}. For \textbf{PINNs} method, we approximate the solution by the feed-forward neural network (FNN) with one input layer, one output layer and $4$ hidden layers with $512$ neurons in each layer, unless otherwise specified. The hyperbolic tangent function (Tanh) is chosen as our activation function. For \textbf{Extrapolation} method used in forecasting problem, we consider the following finite difference formulation for estimating the solution at the $t+1$-th time step:
\begin{equation}
    \frac{\hat{u}(t+1,\bx)-u(t,\bx)}{\Delta t} = \frac{u(t,\bx)-u(t-1,\bx)}{\Delta t}.
\end{equation}
So we have the estimated solution $\hat{u}(t+1, \bx)$ as follows:
\begin{equation}
    \hat{u}(t+1,\bx) = 2u(t,\bx)-u(t-1,\bx).
\end{equation}

\textbf{Training Settings.} 
In general, we train all models with fixed random seed using the L-BFGS optimizer with Strong Wolfe linear search for $500$ iterations and use full batch for most of the following experiments in the numerical experiments unless otherwise specified. All the hyper-parameters are chosen by trial and error.  

\textbf{Loss Design.} For most of our experiments, we consider the spatial and temporal domains to be $[0,2\pi]$ and $[0,1]$ respectively. We choose the collocation points $\{(t_i, x_i)\}$ for $u(t,x)$ in the following way. For spatial points $x_i$, we select the interior points evenly spaced in $[0, 2\pi]$.  For temporal points $t_i$, we select $20$ interior points evenly spaced in the range $[0,1]$. We use the tensor product grid for the collocation points. For simplicity, the weights parameters in equation \eqref{loss_physics} are set to $(w_1, w_2, w_3) = (1,1,1)$. The penalty parameters in equation \eqref{loss_total} are set to $(\lambda_1, \lambda_2) = (1,1)$. 

\textbf{Testing Settings.}
In terms of evaluation for two main tasks: the forward problem, and the forecasting problem, we adopted commonly used metrics in related works \cite{jin2023asymptotic, zhao2023pinnsformer}, including the relative $\ell^2$ or $\ell^{\infty}$ error between the solution approximations $u^{\text{NN}}(t,x,\beta)$ and reference solutions $u^{\text{ref}}(t,x_j,\beta)$, with the relative $\ell^2$ error defined by: 
\begin{equation}
\mathcal{E}(t):=\sqrt{\frac{\sum_j|u^{\text{NN}}(t,x_j,\beta)-u^{\text{ref}}(t,x_j,\beta)|^2}{\sum_j|u^{\text{ref}}(t,x_j,\beta)|^2}}.
\end{equation}

\textbf{Reproducibility.} All models are implemented in PyTorch, and we run our experiments on a server with Intel(R) Xeon(R) Gold 6230 and two A40 48GB GPUs.

\subsection{Linear Convection Equation}
The one-dimensional convection equation is a time-dependent equation widely employed to model transport phenomena. Consider the following general form for the 1D case with initial value condition and periodic boundary condition: 
\begin{equation}\left\{
\begin{aligned}
&\frac{\partial u}{\partial t} + \frac{\partial f(u)}{\partial x} = 0 ,\, (t,x)\in[0,1]\times [0,2\pi],\\
& u(x,0)=\sin(x), \:\:\:  u(0,t)=u(2\pi,t).
\end{aligned} \right.  
\end{equation}
Where $f(u)$ is the flux, $u(t,x)$ is the conserved quantity. Here we set 
$\frac{\mathrm{d} f}{\mathrm{d} u}=50$ be the constant. 
The reference solution is obtained by using a different finite method. For the physics system, we sample $101$ temporal grids and $101$ spatial grids from the above spatial and temporal space similar to \cite{zhao2023pinnsformer}. For the data system, we then use the HSG module to generate extra $5$ temporal stamps (101 spatial points are selected for each temporal stamp) from the grid space for the data system. We use the previous $100$ time steps combined with extra $5$ temporal stamps from the data system for training in forward problems and forecasting problems and use the last time steps for testing in forecasting problems.

\begin{figure}[h]
    \begin{minipage}[b]{0.49\textwidth}
        \centering
        \includegraphics[width=\textwidth]{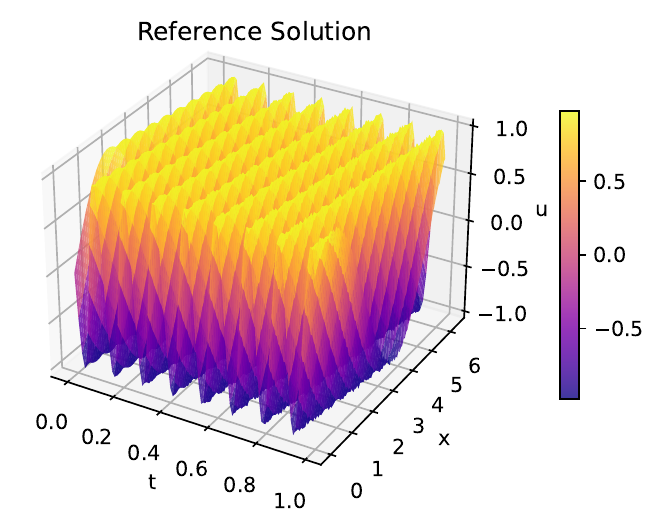} 
    \end{minipage}
    \begin{minipage}[b]{0.49\textwidth}
        \centering
        \includegraphics[width=\textwidth]{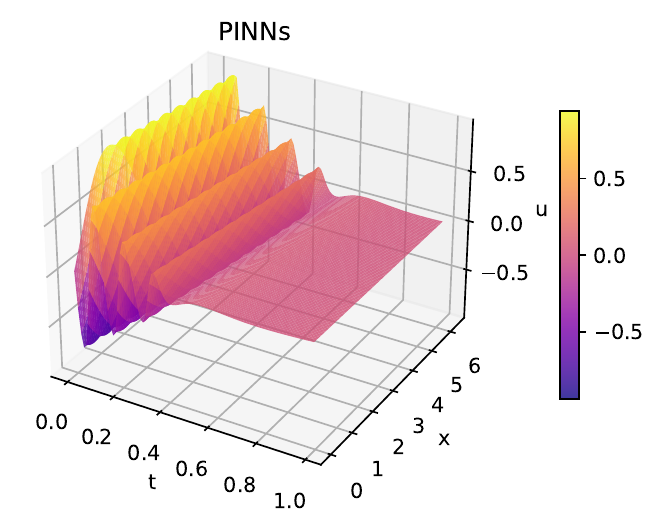} 
    \end{minipage}
    \begin{minipage}[b]{0.49\textwidth}
        \centering
        \includegraphics[width=\textwidth]{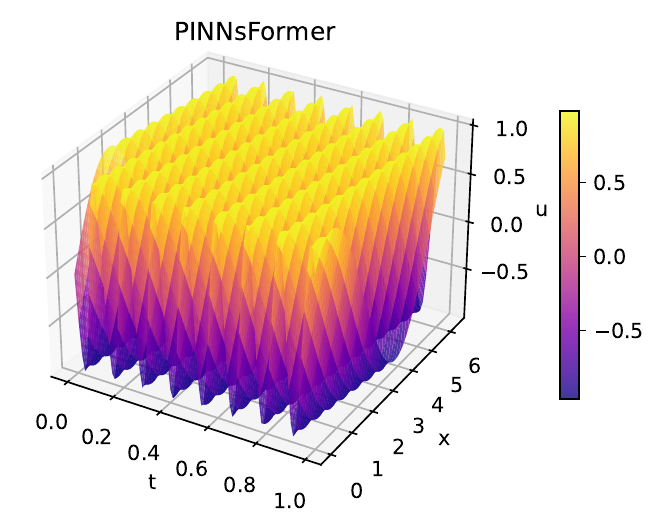} 
    \end{minipage}
    \begin{minipage}[b]{0.49\textwidth}
        \centering
        \includegraphics[width=\textwidth]{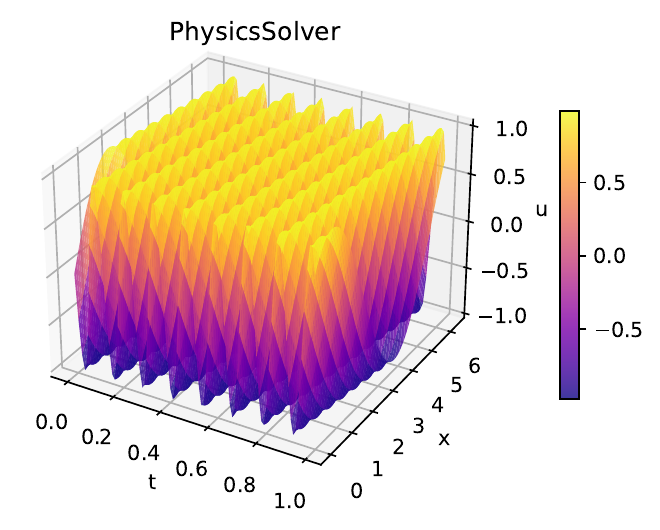} 
    \end{minipage}
    \caption{Performances of methods to the forward problem for the convection equation.}
    \label{fig:test1_forward}
\end{figure}

\begin{figure}[h]
    \centering
    \includegraphics[width=0.32\linewidth]{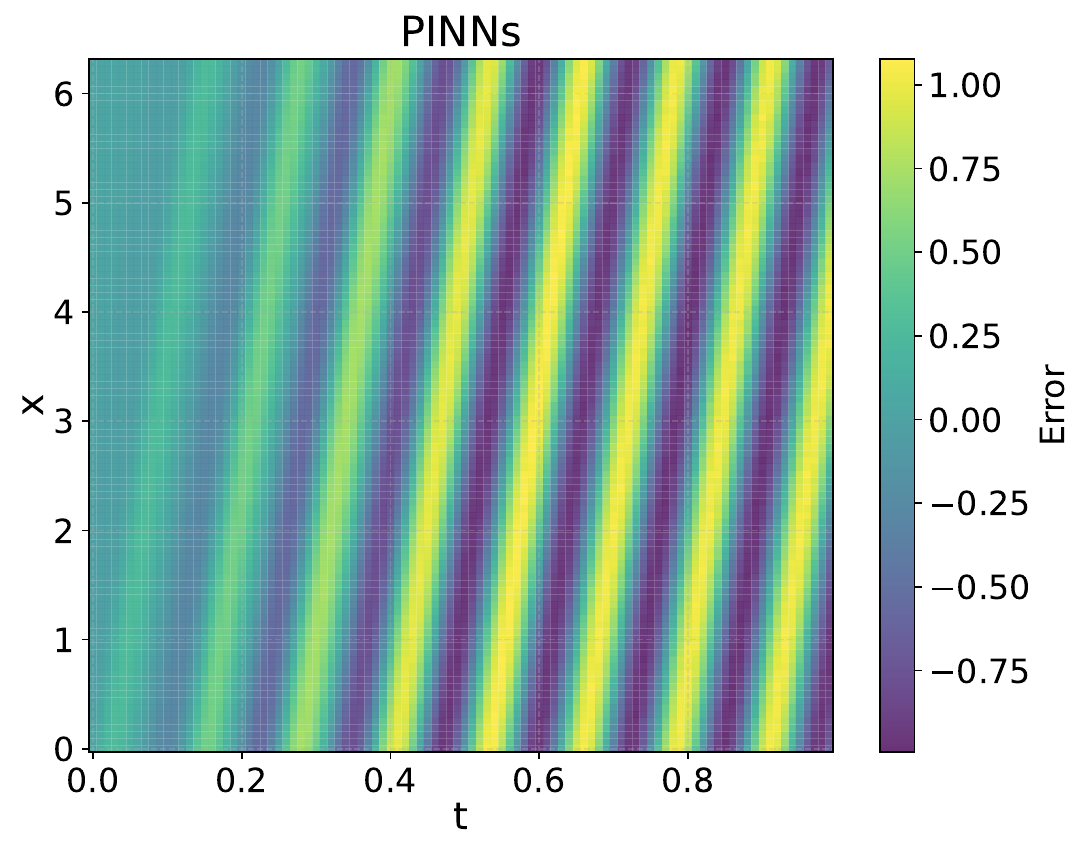}
    \includegraphics[width=0.32\linewidth]{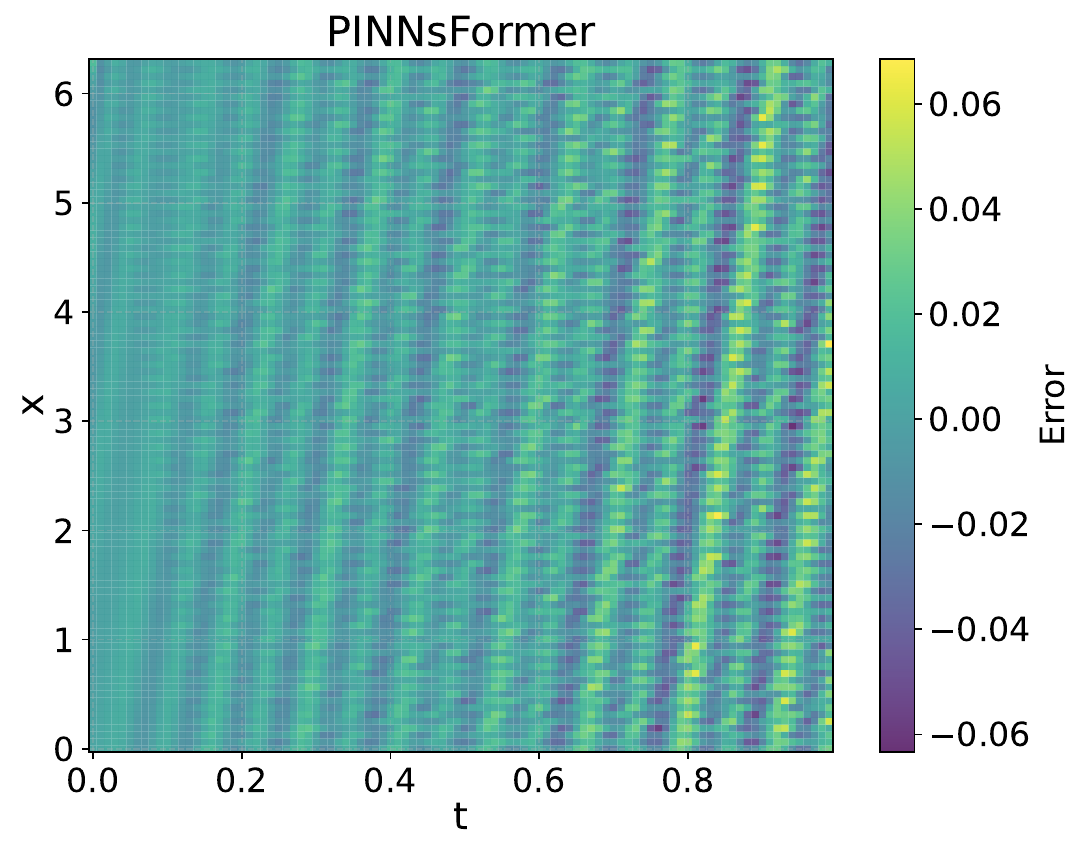}
    \includegraphics[width=0.32\linewidth]{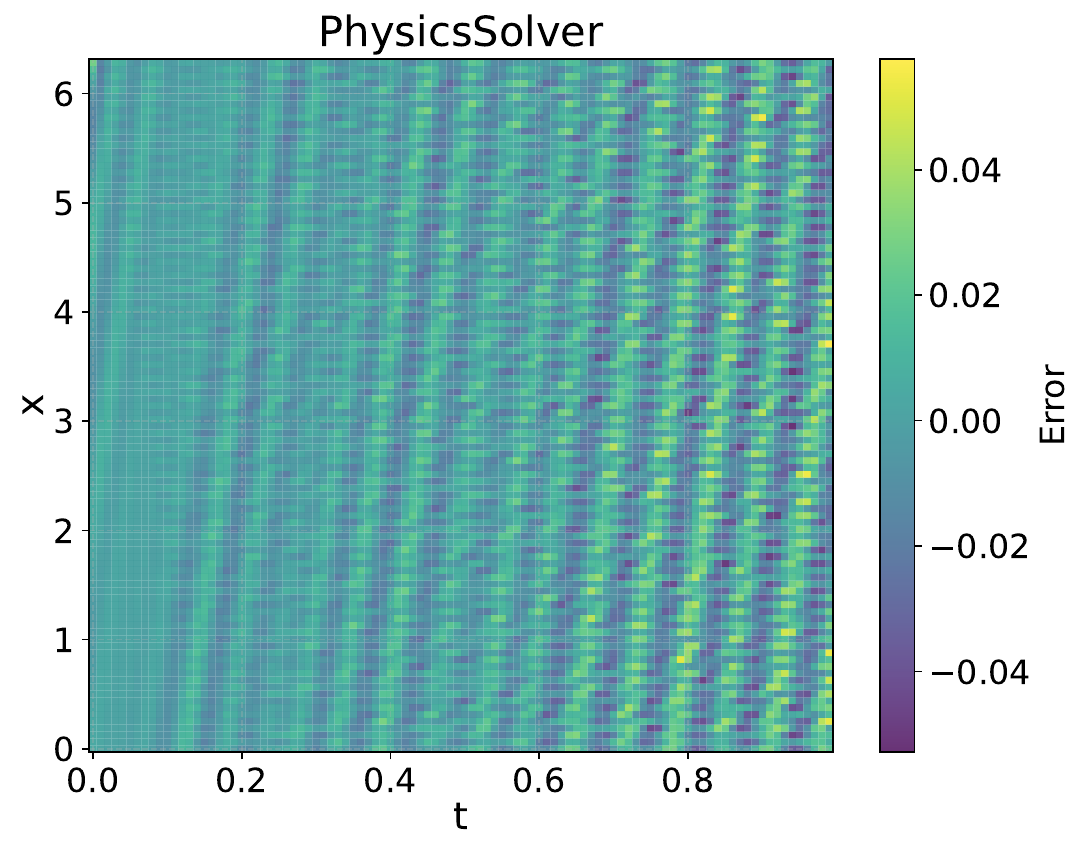}
    \caption{Error distribution across various methods for the convection equation.}
    \label{fig:error_convection}
\end{figure}

\textbf{Forward Problem.}
In the forward problem, we investigate the performance of different methods for the inference of the solution. Figure \ref{fig:test1_forward} shows the comparison for forecasting the solution between different methods. We can find that PINNs fails in capturing the solution with $t$ increases, while PINNsFormer and PhysicsSolver can capture the global solution well. \blue{Furthermore, we show the error distribution across various methods for the convection equation in Figure \ref{fig:error_convection}.} More details can be found in Table \ref{tab:test1_forward}, which shows that PhysicsSolver performs relatively better than PINNsFormer in the forward problem.

\begin{table}[H]
\vspace{-15pt}

	\vskip 0.1in
	\centering
	\begin{small}
		\begin{sc}
			\renewcommand{\multirowsetup}{\centering}
			\setlength{\tabcolsep}{5.5pt}
			\scalebox{1}{
			\begin{tabular}{l|c|c|c}
				\toprule
			    Method  & PINNs & PINNsFormer &PhysicsSolver \\
			    \midrule
                 Error & $8.189 \times 10^{-1}$  & $2.240 \times 10^{-2}$ & $\mathbf{1.558 \times 10^{-2}}$    \\
            \bottomrule
			\end{tabular}}
		\end{sc}
	\end{small}
    \caption{Relative $l^2$ errors for the forward problem.}\label{tab:test1_forward}
\end{table}

\begin{figure}[H]
\centering
\hspace*{-0.9cm}  
\includegraphics[width=1\linewidth]{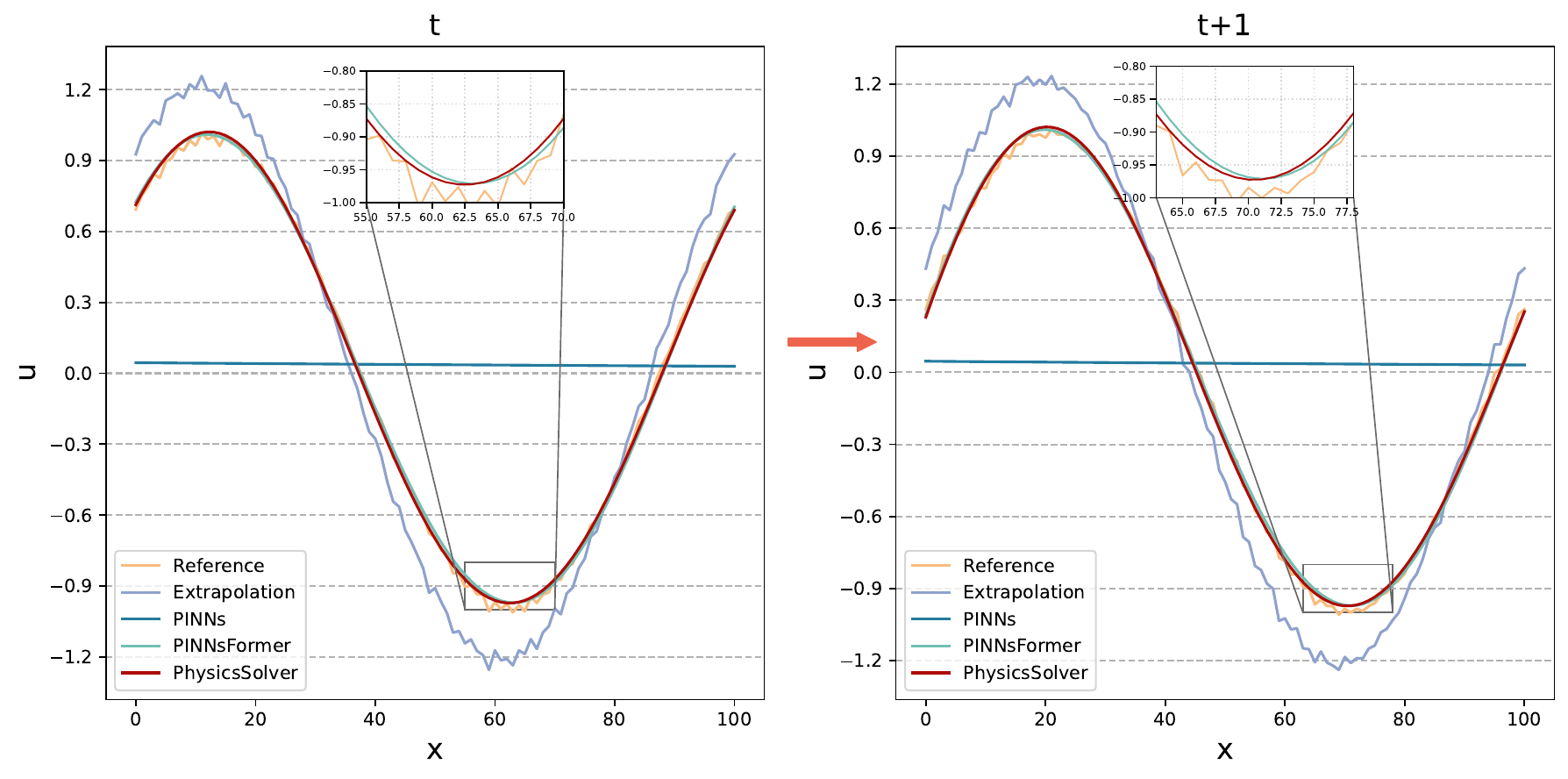}
\caption{The forecasting problem for the convection equation.}
\label{fig:forecast_test1}
\end{figure}

\textbf{Single-step Forecasting Problem.}
In the single-step forecasting problem, we aim to forecast the solution at the next time step. Figure \ref{fig:forecast_test1} compares the solutions obtained by different methods at the last and second-to-last time points. It can be concluded that the Extrapolation method fails to accurately predict the future solution when the solution has relatively significant changes over time. More details can be found in Table \ref{tab:test1_forecast}. We can observe that PhysicsSolver outperforms the other three methods, thanks to its ability to learn a better mapping through the data system, which enhances its predictive capability.

\begin{table}[H]
\vspace{-15pt}
    
	\vskip 0.1in
	\centering
	\begin{small}
		\begin{sc}
			\renewcommand{\multirowsetup}{\centering}
			\setlength{\tabcolsep}{5.5pt}
			\scalebox{1}{
			\begin{tabular}{l|c|c|c|c}
				\toprule
			    Method & Extrapolation & PINNs & PINNsFormer &PhysicsSolver \\
			    \midrule
                 Error &$2.462 \times 10^{-1}$& $9.96 \times 10^{-1}$  & $3.634 \times 10^{-2}$ & $ \mathbf{3.246 \times 10^{-2}}$    \\
            \bottomrule
			\end{tabular}}
		\end{sc}
	\end{small}
    \caption{Relative $l^2$ errors for the forecasting problem.}
	\label{tab:test1_forecast}
\end{table}

\subsection{Reaction Equation}
The one-dimensional reaction equation is an evolution PDE, which is commonly used to model chemical reaction processes. Consider the following  1D equation with  initial value condition and periodic boundary condition: 
\begin{equation}\left\{
\begin{aligned}
    &\frac{\partial u}{\partial t} - \rho u(1-u) = 0,\,(t,x)\in [0,1]\times [0,2\pi],\\
     &u(x,0)=e^{-\frac{(x-\pi)^2}{2(\pi/4)^2}}, \:\:\: u(0,t)=u(2\pi,t).
    \end{aligned}\right.
\end{equation}
where $\rho$ is the reaction coefficient, here we set $\rho=5$ similar to \cite{zhao2023pinnsformer}.

The reference solution is obtained by the following form:
\begin{equation}
    u_{\text{reference}} = \frac{h(x) e^{\rho t}}{h(x)e^{\rho t}+1-h(x)}.
\end{equation}
where $h(x)$ is the function of the initial condition defined in equation \eqref{PDEsys2}.
For the physics system, we sample $101$ temporal grids and $101$ spatial grids from the above spatial and temporal space similar to \cite{zhao2023pinnsformer}. For the data system, we apply the HSG module to generate an extra $10$ temporal stamps (101 spatial points are selected for each temporal stamp) from the grid space. We use the previous $100$ time steps combined with extra $10$ temporal stamps from the data system for training in forward problems and forecasting problems and use the last time steps for testing in forecasting problems.

\begin{figure}[h]
    \begin{minipage}[b]{0.49\textwidth}
        \centering
        \includegraphics[width=\textwidth]{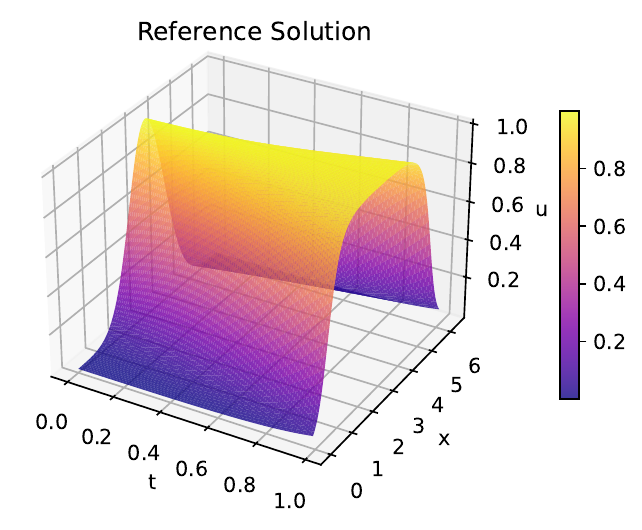} 
    \end{minipage}
    \begin{minipage}[b]{0.49\textwidth}
        \centering
        \includegraphics[width=\textwidth]{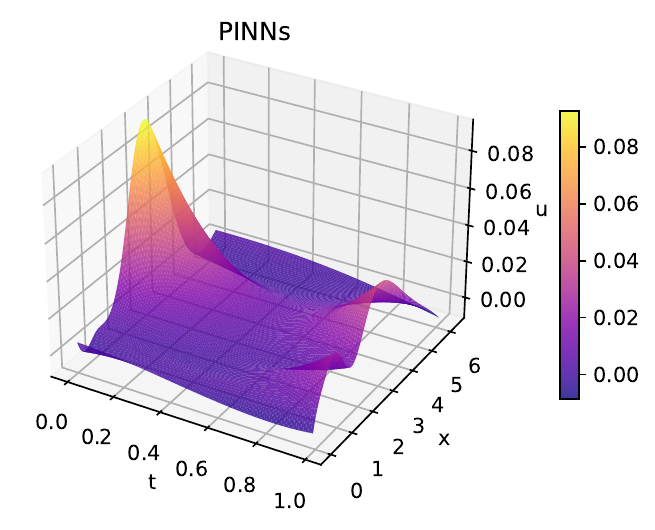} 
    \end{minipage}
    \begin{minipage}[b]{0.49\textwidth}
        \centering
        \includegraphics[width=\textwidth]{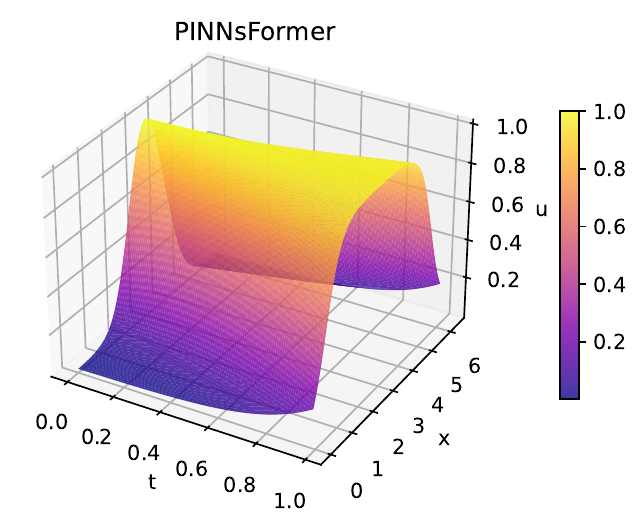} 
    \end{minipage}
    \begin{minipage}[b]{0.49\textwidth}
        \centering
        \includegraphics[width=\textwidth]{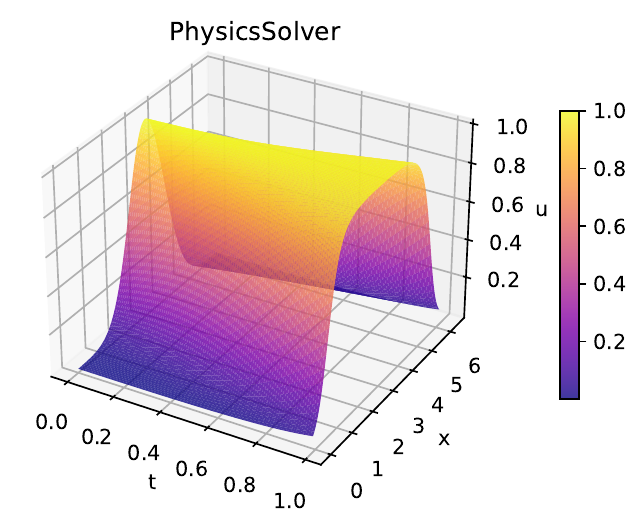} 
    \end{minipage}
    \caption{Performances of methods to the forward problem for the reaction equation.}
    \label{fig:test2_forward}
\end{figure}

\begin{figure}[h]
    \centering
    \includegraphics[width=0.32\linewidth]{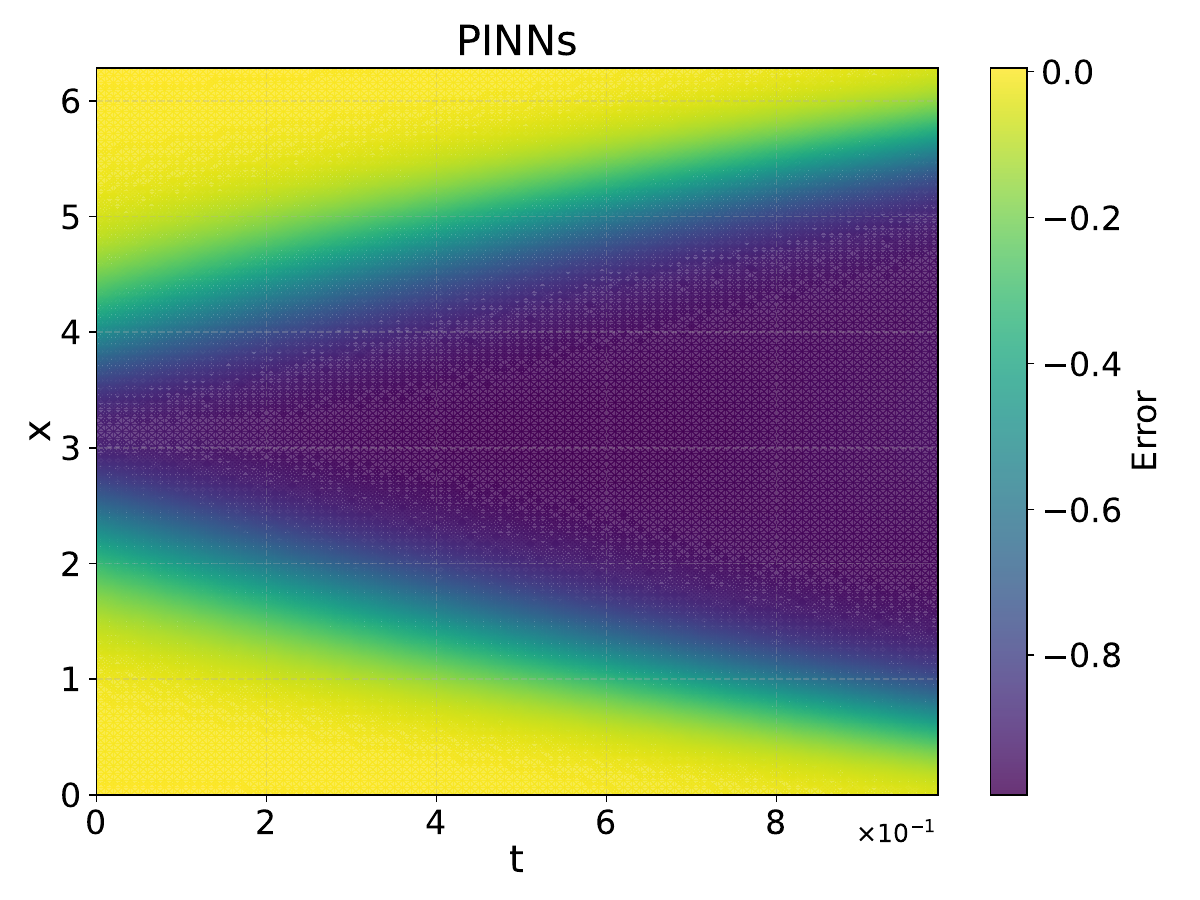}
    \includegraphics[width=0.32\linewidth]{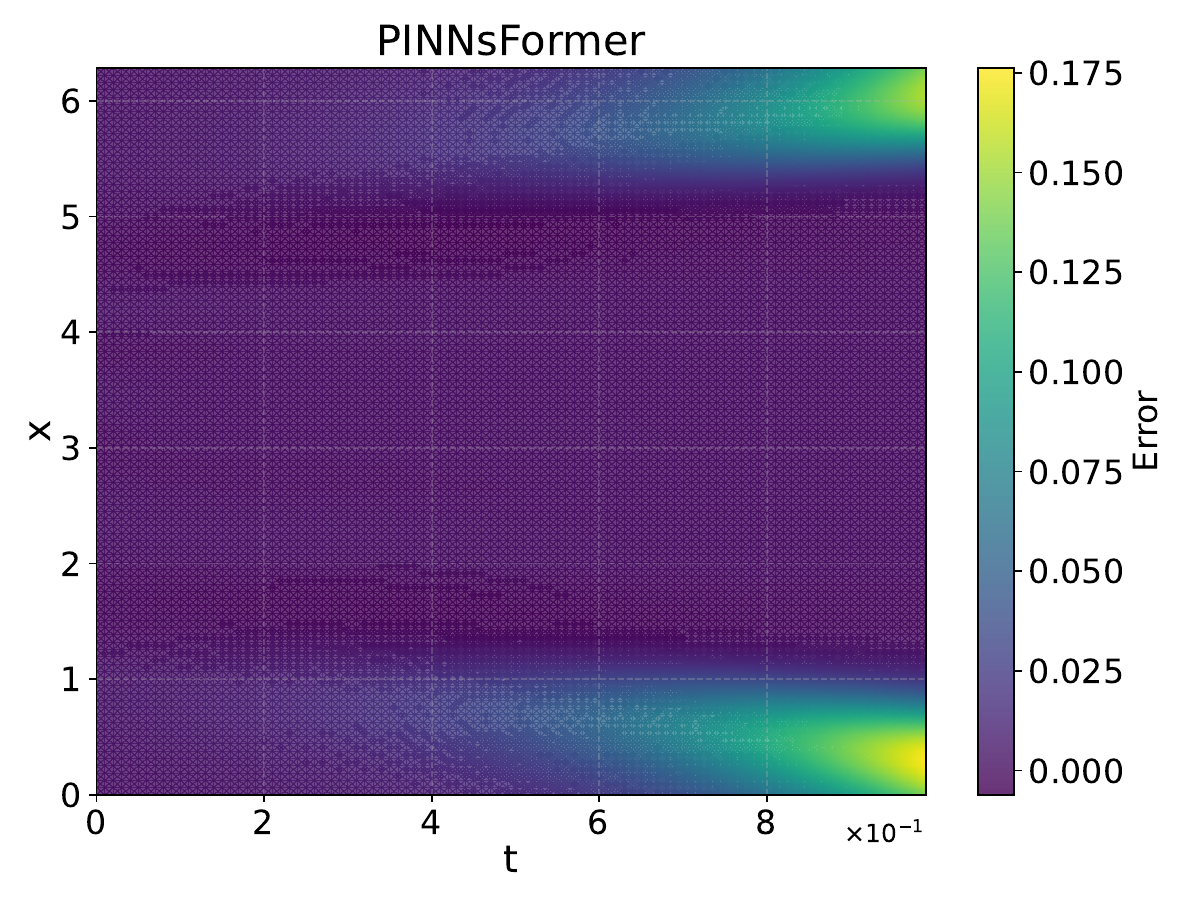}
    \includegraphics[width=0.32\linewidth]{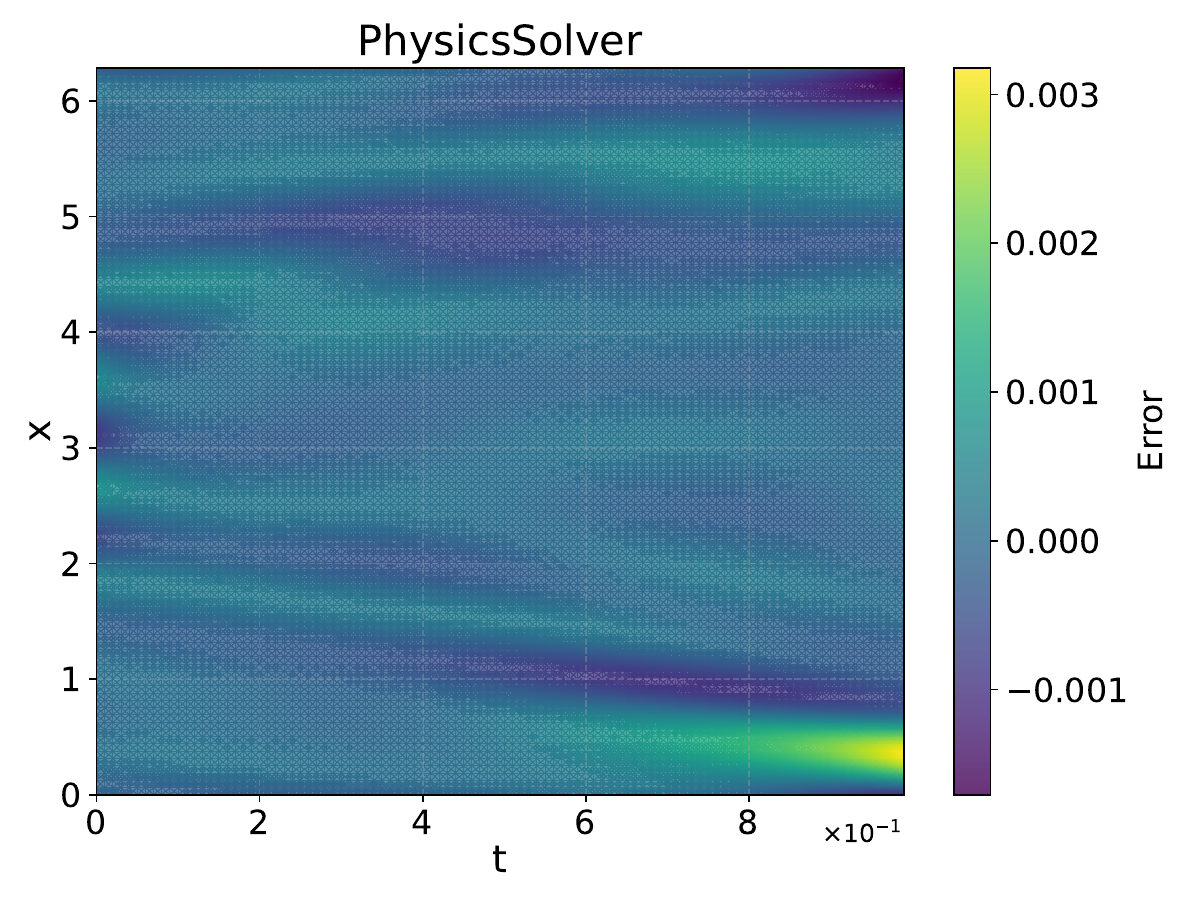}
    \caption{Error distribution across various methods for the reaction equation.}
    \label{fig:error_reaction}
\end{figure}

\textbf{Forward Problem.}
In the forward problem, we investigate the performance of different methods for the inference of the solution. Figure \ref{fig:test2_forward} shows the comparison for forecasting the solution between different methods. We can find that PINNs fails to capture the solution with $t$ increases, while PINNsFormer and PhysicsSolver can capture the global solution well.  \blue{Furthermore, we show the error distribution across various methods for the reaction equation in Figure \ref{fig:error_reaction}.} More details can be found in Table \ref{tab:test2_forward}, it shows that PhysicsSolver performs much better than PINNsFormer in the forward problem thanks to its ability to learn a better mapping through the data system, which can improve the accuracy in the forward problem.\\

\begin{table}[H]
\vspace{-15pt}
    
	\vskip 0.1in
	\centering
	\begin{small}
		\begin{sc}
			\renewcommand{\multirowsetup}{\centering}
			\setlength{\tabcolsep}{5.5pt}
			\scalebox{1}{
			\begin{tabular}{l|c|c|c}
				\toprule
			    Method  & PINNs & PINNsFormer &PhysicsSolver \\
			    \midrule
                 Error & $9.803 \times 10^{-1}$  & $4.563 \times 10^{-2}$ & $\mathbf{5.550 \times 10^{-4}}$    \\
            \bottomrule
			\end{tabular}}
		\end{sc}
	\end{small}
    \caption{Relative $l^2$ errors for the forward problem.}
	\label{tab:test2_forward}
\end{table}

\begin{figure}[h]
\centering
\hspace*{-0.9cm} 
\includegraphics[width=1\linewidth]{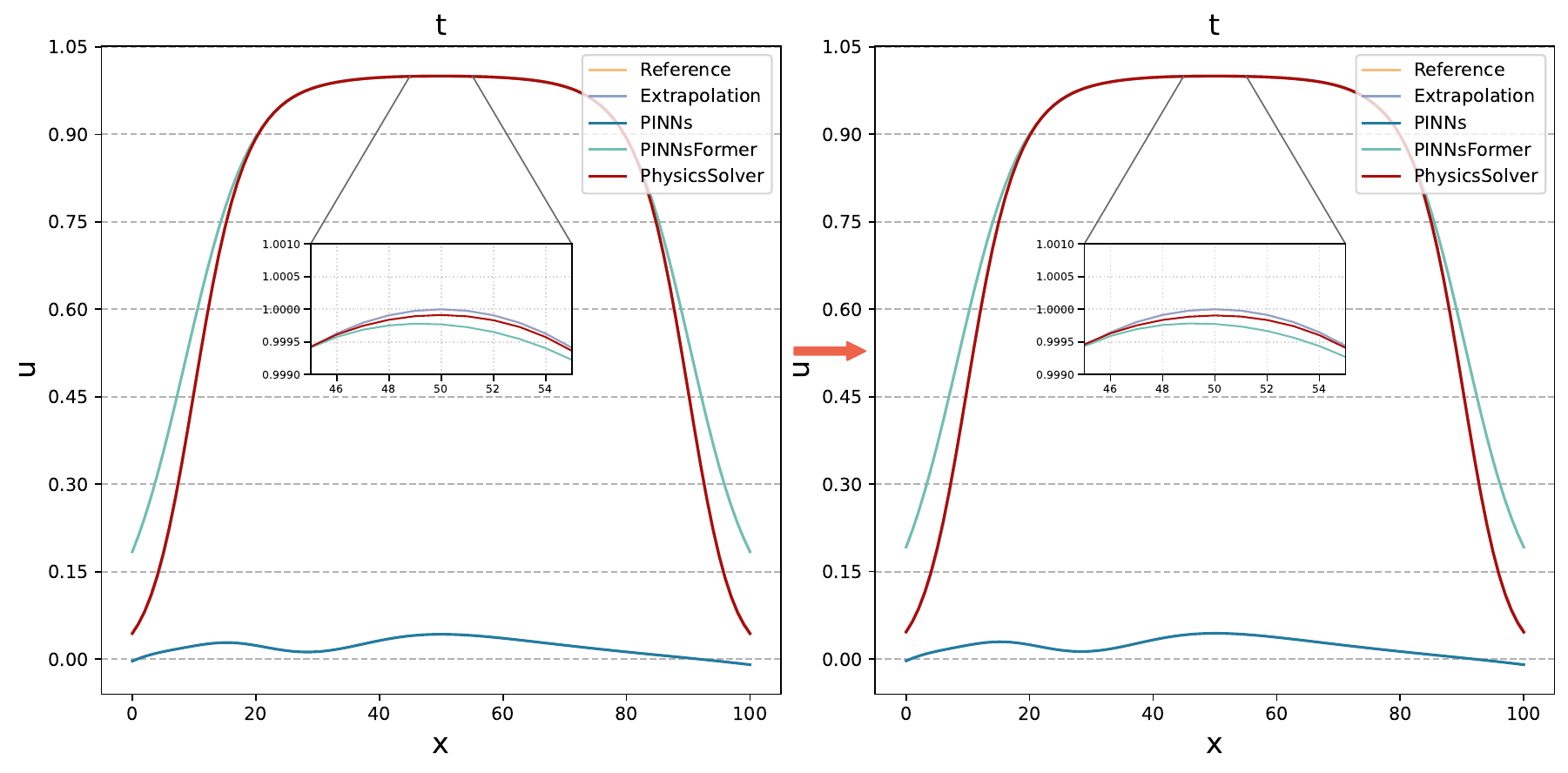}
    \caption{The forecasting problem for the reaction equation.}
    \label{fig:test2_forecast}
\end{figure}

\textbf{Single-step Forecasting Problem.}
In the single-step forecasting problem, we aim to forecast the solution at the next time step. Figure \ref{fig:test2_forecast} compares the solutions obtained by different methods at the last and second-to-last time points. It can be concluded that the PINNs method fails to predict the future solution accurately. However, the Extrapolation method works well when the solution has relatively small changes over time. More details can be found in Table  \ref{tab:test2_forecast}. We can observe that PhysicsSolver outperforms the other three methods. 

\begin{table}[H]
\vspace{-15pt}
    
	\vskip 0.1in
	\centering
	\begin{small}
		\begin{sc}
			\renewcommand{\multirowsetup}{\centering}
			\setlength{\tabcolsep}{5.5pt}
			\scalebox{1}{
			\begin{tabular}{l|c|c|c|c}
				\toprule
			    Method & Extrapolation & PINNs & PINNsFormer &PhysicsSolver \\
			    \midrule
                 Error &$6.915 \times 10^{-3}$& $9.726 \times 10^{-1}$  & $8.590 \times 10^{-2}$ & $\mathbf{9.460 \times 10^{-4}}$    \\
            \bottomrule
			\end{tabular}}
		\end{sc}
	\end{small}
    \caption{Relative $l^2$ errors for the forecasting problem.}
	\label{tab:test2_forecast}
\end{table}

\subsection{Heat Diffusion Equation}
The one-dimensional heat equation is the prototypical parabolic partial differential equation to model how a quantity, such as heat, diffuses across a given region. Consider the following general form for the 1D case with periodic boundary conditions: 
\begin{equation}
\left\{\begin{aligned}
    &\frac{\partial u}{\partial t} = \alpha \frac{\partial^2 u}{\partial x^2},\,(t,x)\in [0,0.2]\times[0,1],\\
    &u(x,0)=\sin(\pi x), \:\:\:  u(0,t) = u(1,t) = 0.
    \end{aligned}\right.
\end{equation}
where $\alpha$ is a positive coefficient called the thermal diffusivity of the medium, here we set $\alpha = 1$ for simplicity.

The reference solution is obtained by the following form:
\begin{equation}
    u_{\text{reference}} =  e^{-\pi^2 t}\sin(\pi x).
\end{equation}
For the physics system, we sample $5$ temporal grids (or temporal stamps) and $101$ spatial grids from the above spatial and temporal space. For the data system, we use the HSG module to generate extra $2$ temporal stamps ($101$ spatial points are selected for each temporal stamp) from the grid space for the data system. We use the previous $4$ time steps combined with extra 2 temporal stamps from the data system for training in forward problems and forecasting problems and use the last time steps for testing in forecasting problems. 

\begin{figure}[h]
\centering
    \begin{minipage}[b]{0.49\textwidth}
        \centering
        \includegraphics[width=\textwidth]{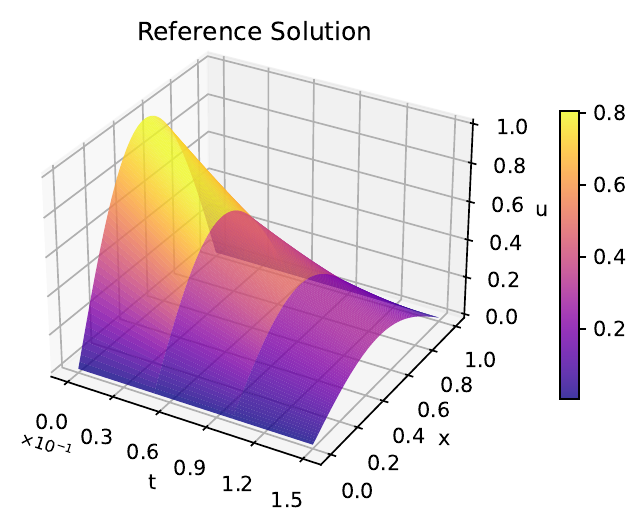} 
    \end{minipage}
    \begin{minipage}[b]{0.49\textwidth}
        \centering
        \includegraphics[width=\textwidth]{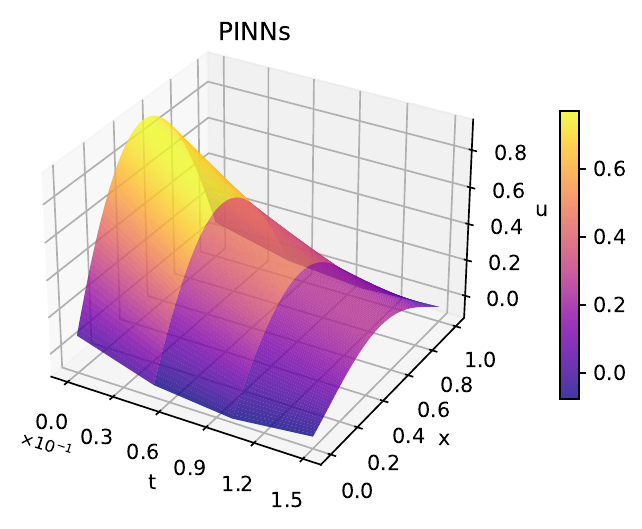} 
    \end{minipage}
    \begin{minipage}[b]{0.49\textwidth}
        \centering
        \includegraphics[width=\textwidth]{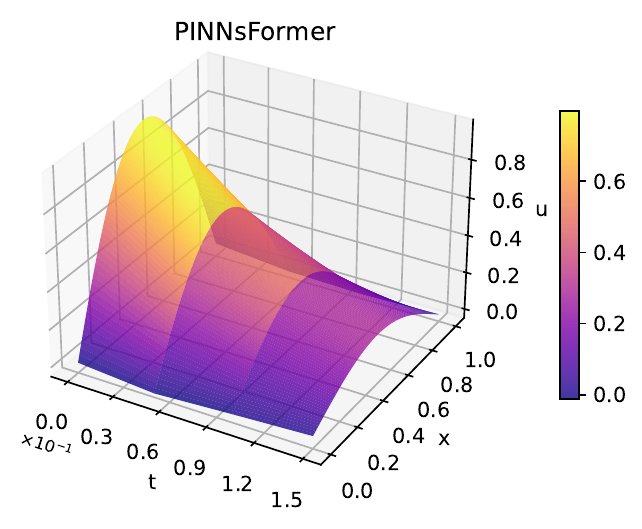} 
    \end{minipage}
    \begin{minipage}[b]{0.49\textwidth}
        \centering
        \includegraphics[width=\textwidth]{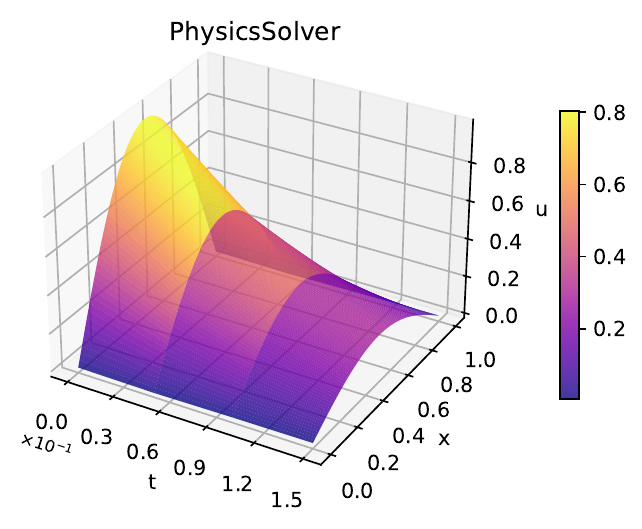} 
    \end{minipage}
    \caption{Performances of methods to the forward problem for the heat equation.}
    \label{fig:test3_forward}
\end{figure}

\begin{figure}[h]
    \centering
    \includegraphics[width=0.32\linewidth]{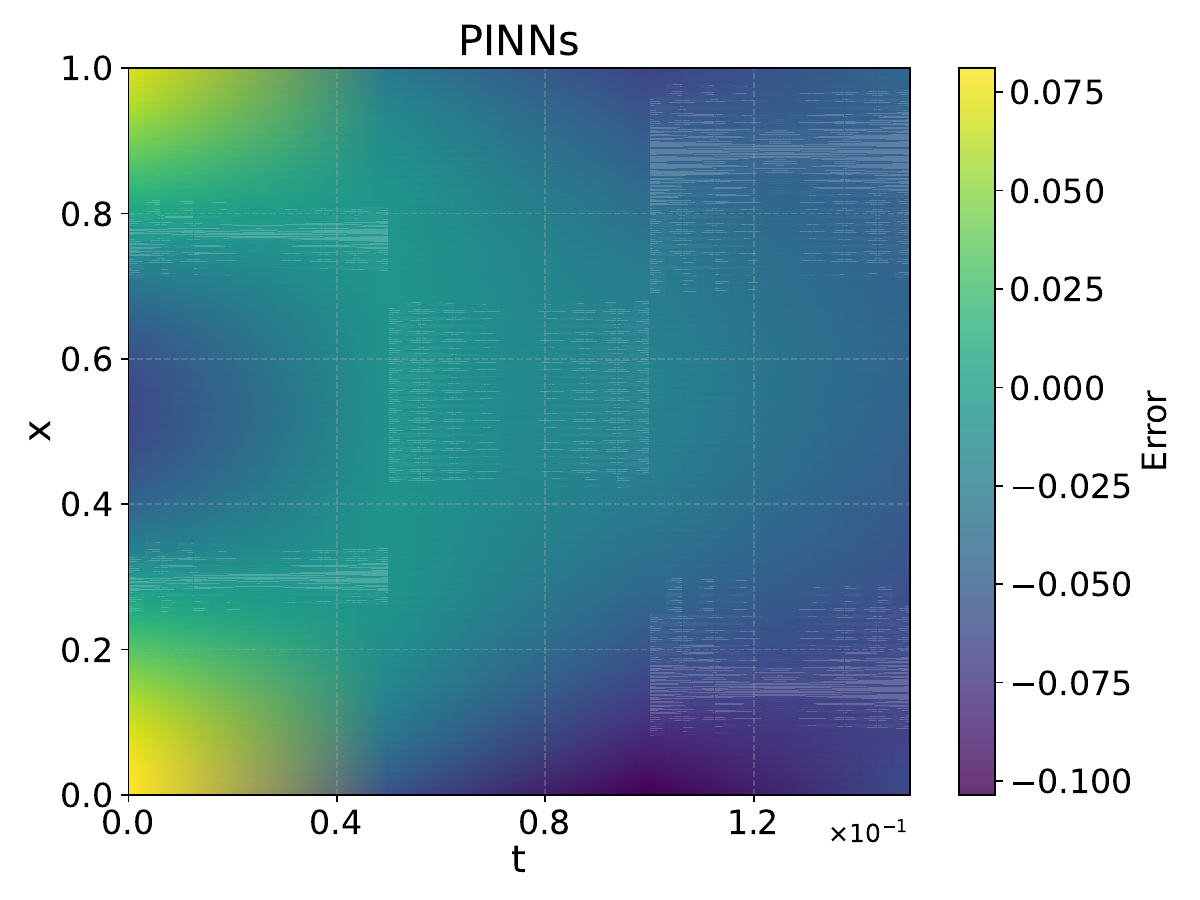}
    \includegraphics[width=0.32\linewidth]{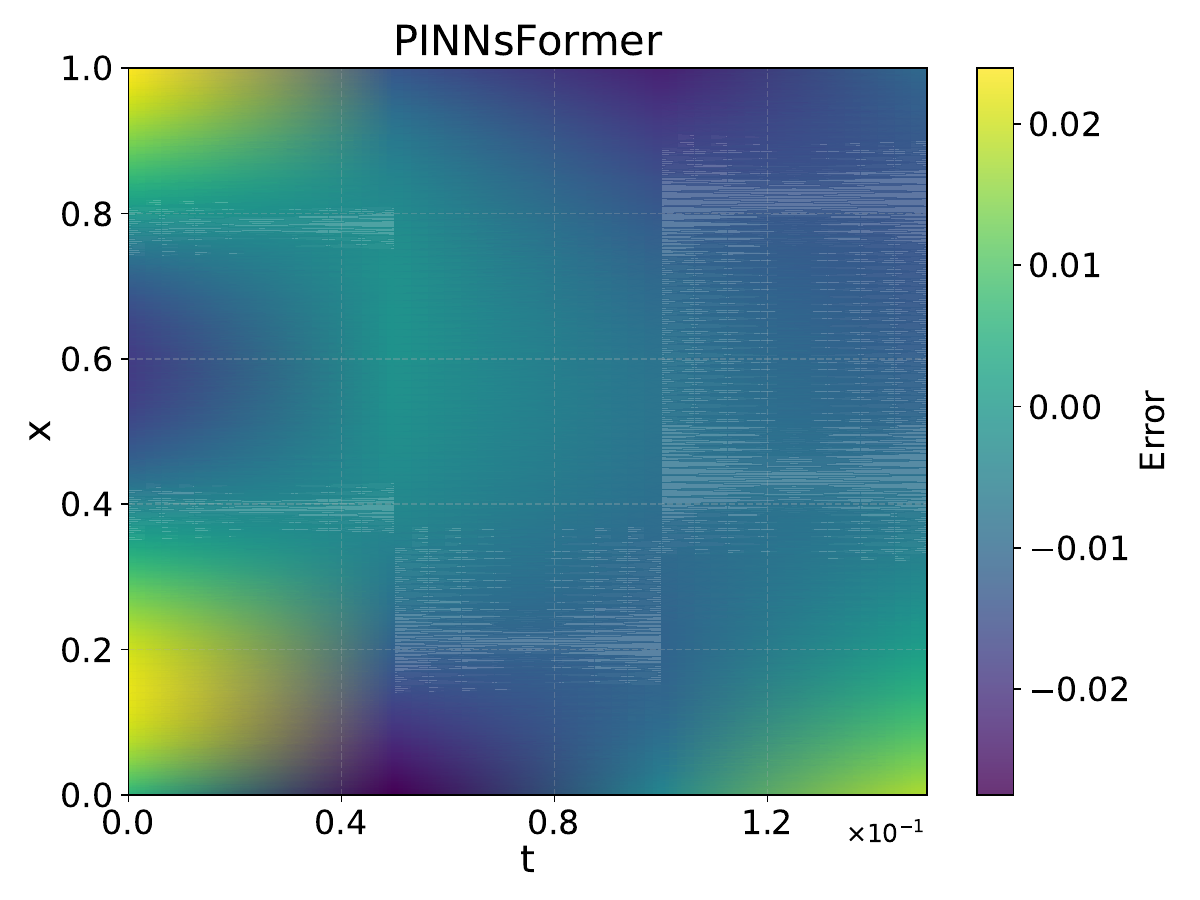}
    \includegraphics[width=0.32\linewidth]{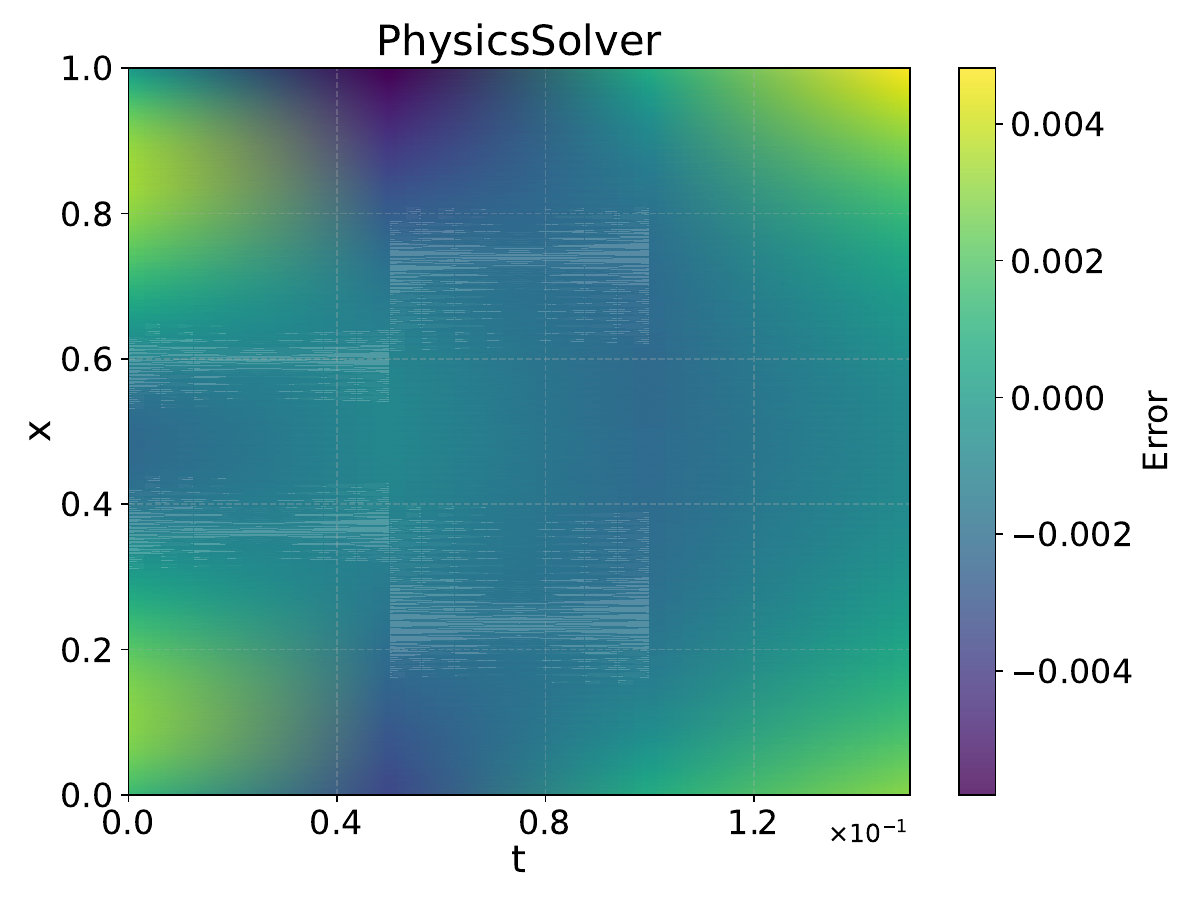}
    \caption{Error distribution across various methods for the heat equation.}
    \label{fig:error_heat}
\end{figure}

\textbf{Forward Problem.}
In the forward problem, we investigate the performance of different methods for the inference of the solution. Figure \ref{fig:test3_forward} shows the comparison for forecasting the solution between different methods. We can find that both the three methods can capture the solution as $t$ increases, while PINNsFormer and PhysicsSolver can capture the global solution better. \blue{Furthermore, we show the error distribution across various methods for the heat equation in Figure \ref{fig:error_heat}.} More details can be found in Table \ref{tab:test3_forward}, which shows that PhysicsSolver performs much better than PINNsFormer in the forward problem thanks to its ability to learn a better mapping through the data system, which improves the accuracy in the forward problem.

\begin{table}[H]
\vspace{-15pt}
    
	\vskip 0.1in
	\centering
	\begin{small}
		\begin{sc}
			\renewcommand{\multirowsetup}{\centering}
			\setlength{\tabcolsep}{5.5pt}
			\scalebox{1}{
			\begin{tabular}{l|c|c|c}
				\toprule
			    Method  & PINNs & PINNsFormer &PhysicsSolver \\
			    \midrule
                 Error & $1.000 \times 10^{-1}$  & $2.702 \times 10^{-2}$ & $\mathbf{4.622 \times 10^{-3}}$    \\
            \bottomrule
			\end{tabular}}
		\end{sc}
	\end{small}
    \caption{Relative $l^2$ errors for the forward problem.}
	\label{tab:test3_forward}
\end{table}

\begin{figure}[h]
\centering
\hspace{-0.9cm}\includegraphics[width=1\linewidth]{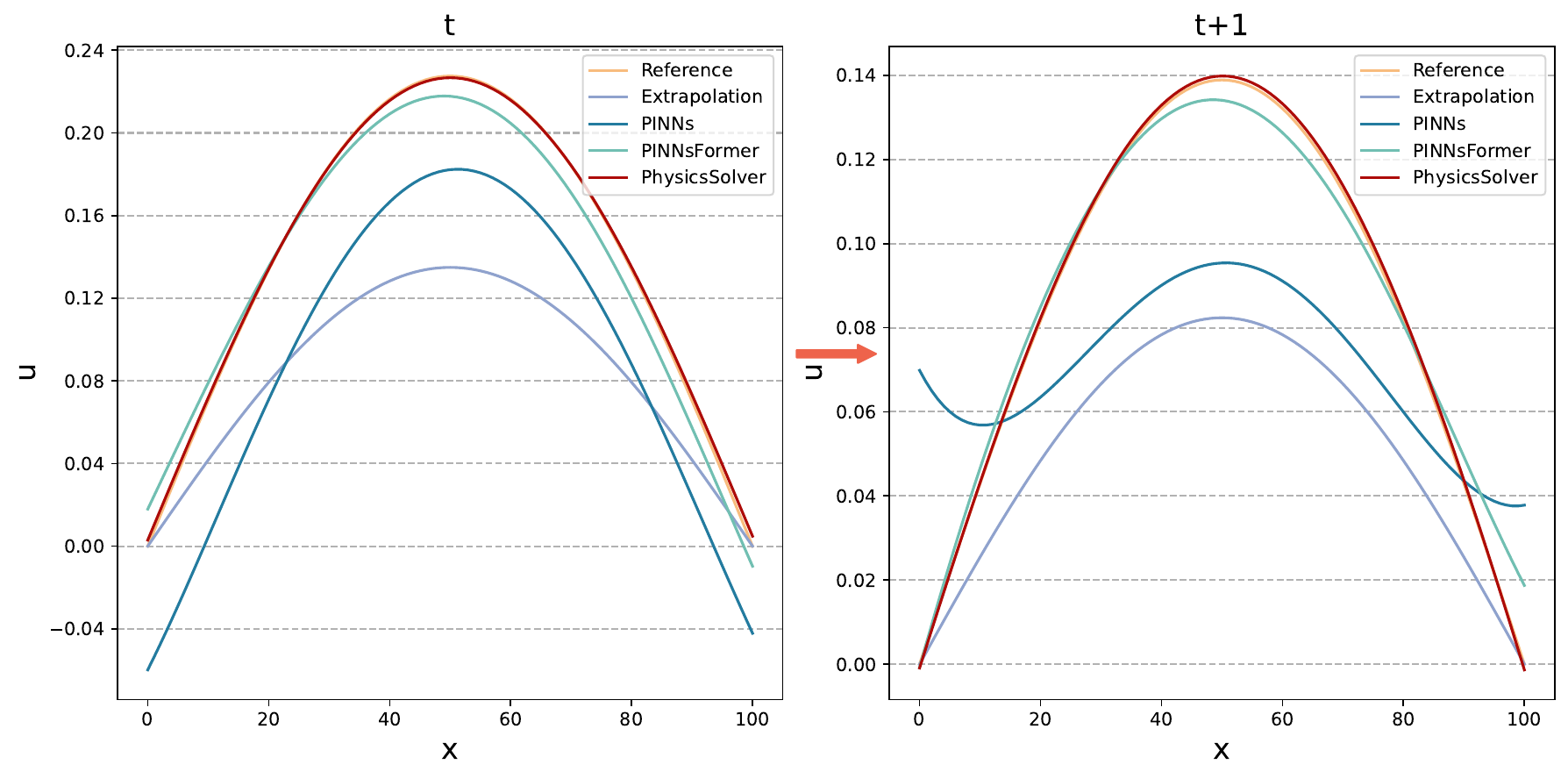}
    \caption{Forecasting problem for the heat equation.}
    \label{fig:test3_forecast}
\end{figure}

\textbf{Single-step Forecasting Problem.}
Similarly, we aim to forecast the solution at the next time step. Figure \ref{fig:test3_forecast} compares the solutions obtained by different methods at the last and second-to-last time points. It can be concluded that the PINNs fails to accurately predict the future solution because it can not learn the time-dependent relationship between different solutions. The Extrapolation method fails to accurately predict the future solution when the solution has relatively significant changes over time. Also, PINNsformer can not capture the solution well since the model ignores the available data information. More details can be found in Table  \ref{tab:test3_forecast}. We can observe that PhysicsSolver outperforms the other three methods. 
\begin{table}[H]
\vspace{-15pt}
    
	\vskip 0.1in
	\centering
	\begin{small}
		\begin{sc}
			\renewcommand{\multirowsetup}{\centering}
			\setlength{\tabcolsep}{5.5pt}
			\scalebox{1}{
			\begin{tabular}{l|c|c|c|c}
				\toprule
			    Method & Extrapolation & PINNs & PINNsFormer &PhysicsSolver \\
			    \midrule
                 Error &$6.380 \times 10^{-1}$& $3.241 \times 10^{-1}$  & $6.560 \times 10^{-2}$ & $\mathbf{1.070 \times 10^{-3}}$    \\
            \bottomrule
			\end{tabular}}
		\end{sc}
	\end{small}
    \caption{Relative $l^2$ errors for the forecasting problem.}
	\label{tab:test3_forecast}
\end{table}
\subsection{2D Navier-Stokes Equations}
Navier–Stokes equations mathematically represent the mass and momentum conservation laws for Newtonian fluids that link pressure, velocity, temperature and density. The Navier–Stokes equations, both in their full and simplified forms, play a crucial role in the design of aircraft and automobiles, the study of blood flow, the development of power stations, pollution analysis, and a wide range of other applications. When coupled with Maxwell's equations, they can also be used to model and analyze magnetohydrodynamics. Consider the following system of parabolic PDEs, which models the behaviour of incompressible fluid flow in two-dimensional space.
\begin{equation}\left\{
\begin{aligned}
&u_{x}+v_{y}=0,\\ 
    &u_t + \lambda_1 [(u ,v) \cdot \nabla u]  = - p_x + \lambda_2 \Delta u,  \\
    &v_t + \lambda_1 [(u,v) \cdot \nabla v] = - p_y + \lambda_2 \Delta v.
\end{aligned}\right.
\end{equation}
Where $u(t,x,y)$ is the first component of the velocity field $\boldsymbol{u}$, and $v(t,x,y)$ is the second component of the velocity field $\boldsymbol{u}$, that is, $\boldsymbol{u}=(u,v)$. $p(t,x,y)$ is the pressure. $\nabla u = (\frac{\partial u}{\partial x}$, $\frac{\partial u}{\partial y}), \, \Delta u = \frac{\partial^2 u}{\partial x^2} + \frac{\partial^2 u}{\partial y^2}$, $\nabla v = (\frac{\partial v}{\partial x}$, $\frac{\partial v}{\partial y}), \, \Delta v = \frac{\partial^2 v}{\partial x^2} + \frac{\partial^2 v}{\partial y^2}$. 
$\lambda = (\lambda_1, \lambda_2)$ are the unknown parameters. Here, we set $\lambda_1 = 1$ and $\lambda_2 = 0.01$. The system does not have an analytical solution, we use the simulated solution provided by~\cite{raissi2019physics} to be the reference solution. In practice, we set the spatial space $x \in [1,7]$, $y \in [-2,2]$, the temporal space $t \in [0, 1]$. 

For the physics system, we sample $10$ temporal grids (or temporal stamps) and $250$ spatial grids from the above spatial and temporal space respectively. For the data system, we use the HSG module to generate extra $4$ temporal stamps ($250$ spatial points are selected for each temporal stamp) from the grid space for the data system. We use the previous $9$ time steps combined with extra 2 temporal stamps from the data system for training in the forward and the forecasting problem and use the last time steps for testing in the forecasting problem.

\begin{figure}[h]
\centering
    \begin{minipage}[b]{0.49\textwidth}
        \centering
        \includegraphics[width=\textwidth]{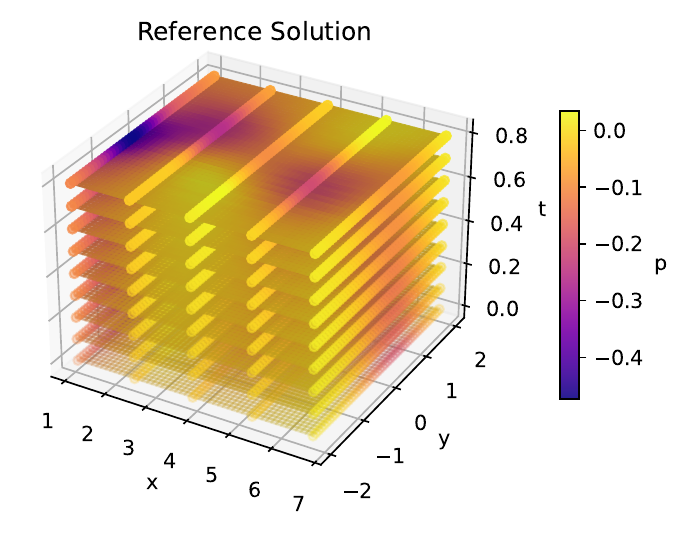} 
    \end{minipage}
    \begin{minipage}[b]{0.49\textwidth}
        \centering
        \includegraphics[width=\textwidth]{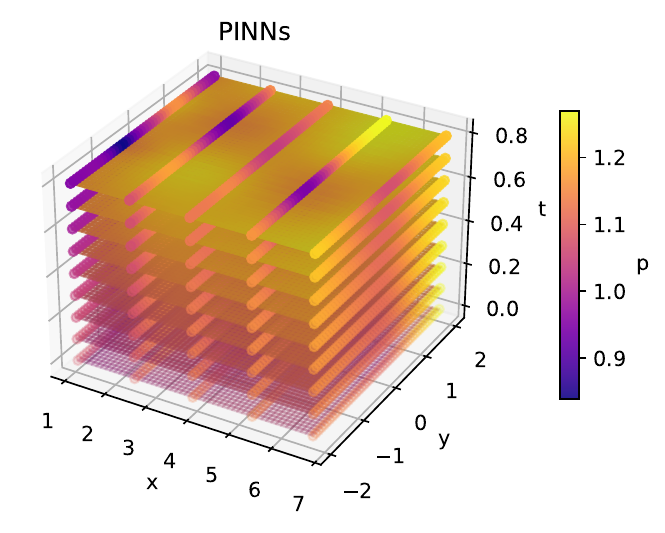} 
    \end{minipage}
    \begin{minipage}[b]{0.49\textwidth}
        \centering
        \includegraphics[width=\textwidth]{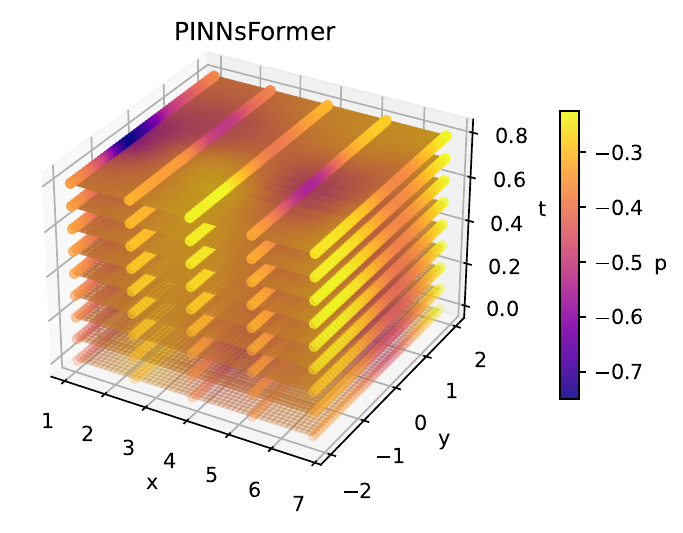} 
    \end{minipage}
    \begin{minipage}[b]{0.49\textwidth}
        \centering
        \includegraphics[width=\textwidth]{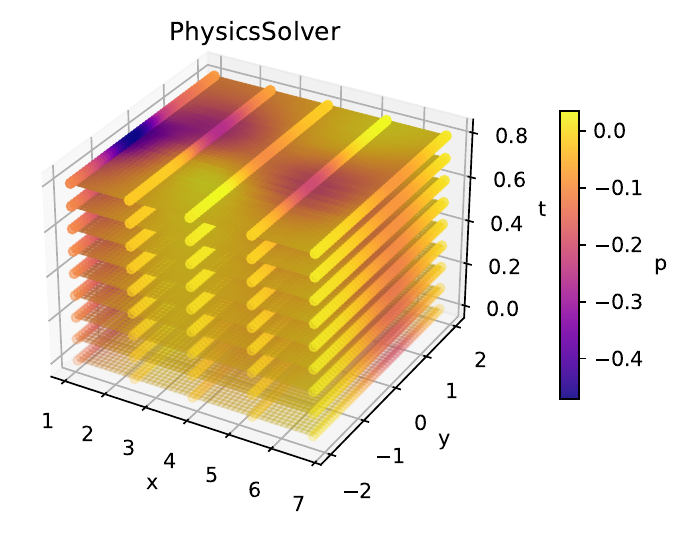} 
    \end{minipage}
    \caption{Forward problem for the Navier-Stokes equation.}
    \label{fig:test4_forward}
\end{figure}

\begin{figure}[h]
    \centering
    \includegraphics[width=0.32\linewidth]{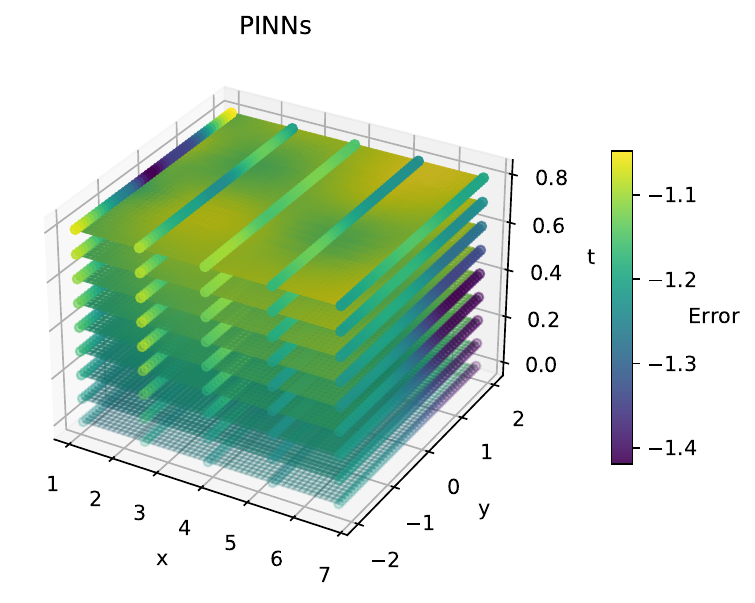}
    \includegraphics[width=0.32\linewidth]{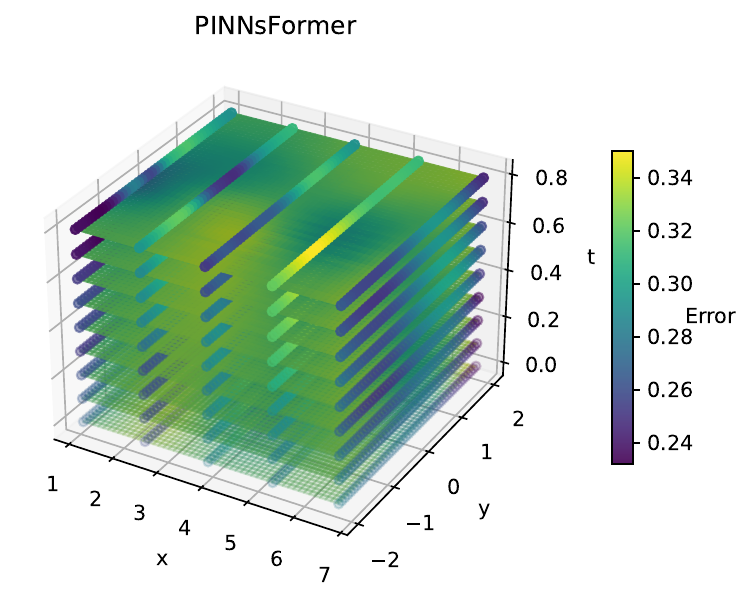}
    \includegraphics[width=0.32\linewidth]{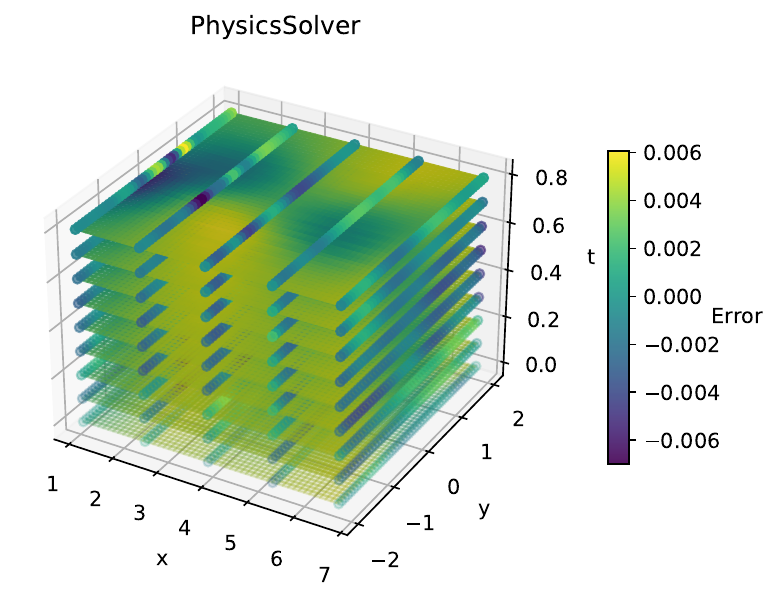}
    \caption{Error distribution across various methods for the Navier-Stokes equation.}
    \label{fig:error_ns}
\end{figure}

\textbf{Forward Problem.}
In the forward problem, we investigate the performance of different methods for the inference of the solution for pressure $p$. Figure \ref{fig:test4_forward} shows the comparison for forecasting the solution between different methods. We can find that both the PINNs method and PINNsFormer method can not capture the solution well at different time stamps, while PhysicsSolver can capture the global solution better. \blue{Furthermore, we show the error distribution across various methods for the Navier-Stokes equation in Figure \ref{fig:error_ns}.} More details can be found in Table \ref{tab:test4_forward}, which shows that PhysicsSolver performs much better than PINNs and PINNsFormer in the forward problem thanks to its ability to learn a better mapping through the data system. 
 
It can also be concluded that the incorporation of the data system loss is more effective for the forward problem when the solution exhibits little change across different time stamps. 

\begin{table}[H]
\vspace{-15pt}
    
	\vskip 0.1in
	\centering
	\begin{small}
		\begin{sc}
			\renewcommand{\multirowsetup}{\centering}
			\setlength{\tabcolsep}{5.5pt}
			\scalebox{1}{
			\begin{tabular}{l|c|c|c}
				\toprule
			    Method  & PINNs & PINNsFormer &PhysicsSolver \\
			    \midrule
                 Error & $5.824 \times 10^{0}$  & $2.139 \times 10^{0}$ & $\mathbf{9.224 \times 10^{-3}}$    \\
            \bottomrule
			\end{tabular}}
		\end{sc}
	\end{small}
    \caption{Relative $l^2$ errors.}
	\label{tab:test4_forward}
\end{table}

\begin{figure}[h]
\centering
\hspace*{-0.9cm}
\includegraphics[width=1\linewidth]{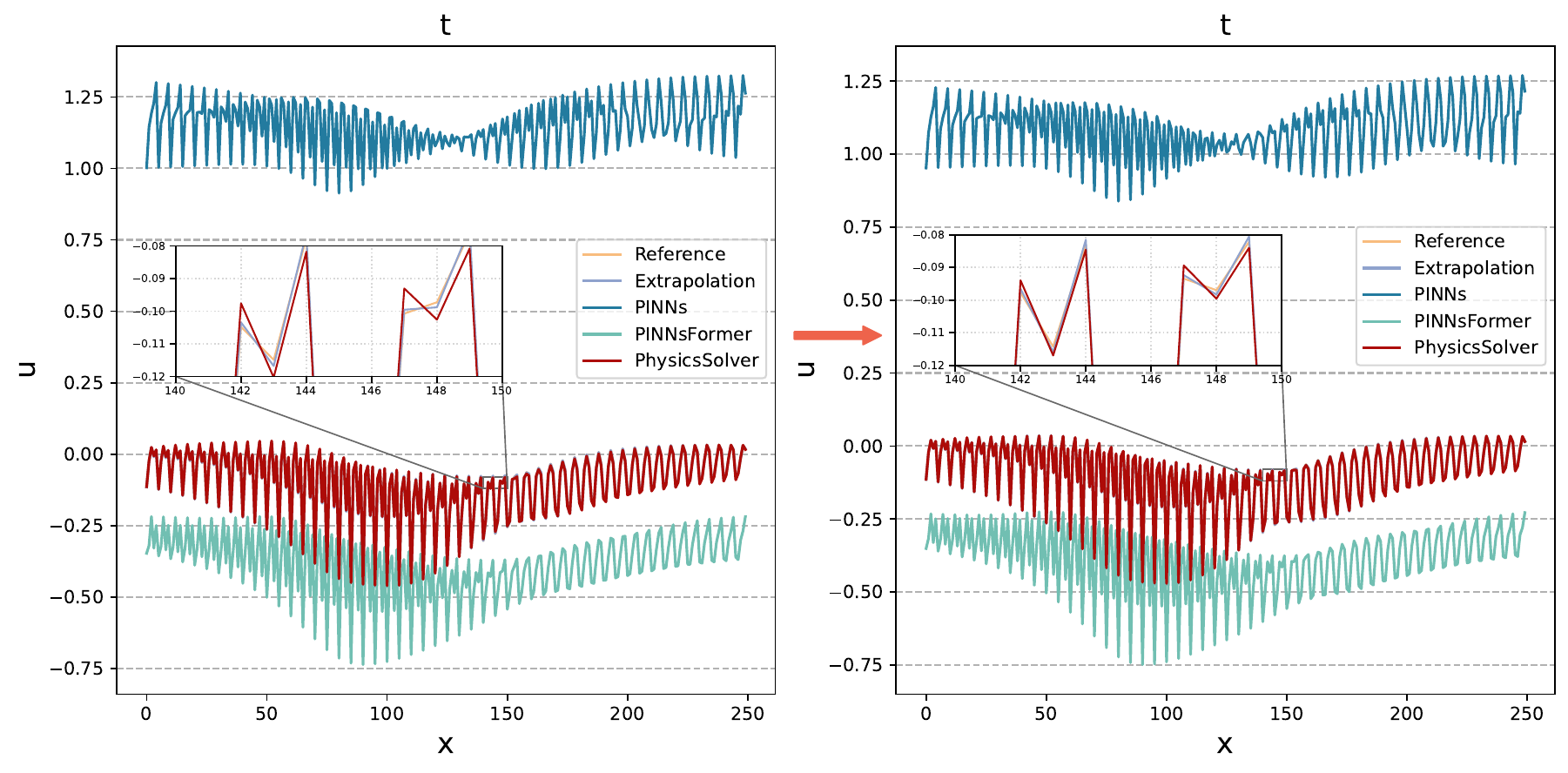}
    \caption{Forecasting problem for the Navier-Stokes equation.}
    \label{fig:test4_forecast}
\end{figure}

\textbf{Single-step Forecasting Problem.}
In the single-step forecasting problem, we also aim to forecast the solution at the next time step. Figure \ref{fig:test4_forecast} compares the solutions obtained by different methods at the last and second-to-last time points.
It can be concluded that the PINNs fails to predict the future solution accurately since it can not learn the time-dependent relationship between different solutions well when the length of time is relatively large. At the same time, PINNsFormer can not capture the solution well since it overlooks the available data information. However, the Extrapolation method can accurately predict the future solution since the solution has relatively small changes over time. More details can be found in Table \ref{tab:test4_forecast}. We can observe that PhysicsSolver outperforms the other two methods and has a similar accuracy to the Extrapolation method. 

\begin{table}[H]
\vspace{-15pt}
    
	\vskip 0.1in
	\centering
	\begin{small}
		\begin{sc}
			\renewcommand{\multirowsetup}{\centering}
			\setlength{\tabcolsep}{5.5pt}
			\scalebox{1}{
			\begin{tabular}{l|c|c|c|c}
				\toprule
			    Method & Extrapolation & PINNs & PINNsFormer &PhysicsSolver \\
			    \midrule
                 Error &$\mathbf{7.437 \times 10^{-3}}$& $7.796 \times 10^{0}$  & $1.912 \times 10^{0}$ & $\mathbf{1.384 \times 10^{-2}}$    \\
            \bottomrule
			\end{tabular}}
		\end{sc}
	\end{small}
    \caption{Relative $l^2$ errors.}
	\label{tab:test4_forecast}
\end{table}

\blue{\subsection{Special Case: multi-step Forecasting Problem for the Reaction Equation}}

\blue{In the single-step forecasting problem, $\tilde{t}$ is a single time step, which is one time step later than previous time $t$ used for learning. In \textit{multi-steps forecasting problem}, suppose the time steps involved to train are $\{t_{m}\}_{m=1}^{m_1}$, and the time steps to be predicted are denoted by $\{t_{m}\}_{m=m_1+1}^{m_2}$. There are several different methods for multi-steps forecasting, such as direct multi-step forecast strategy, recursive multi-step forecast strategy, hybrid forecast strategies, we refer the readers to \cite{taieb2012recursive, taieb2010multiple}. Here we study the direct multi-step forecast strategy for PINNs, PINNsFormer and PhysicsSolver, which means we directly generate the solutions on $\{t_{m}\}_{m=m_1+1}^{m_2}$ time steps by utilizing the 
pre-trained models. For the extrapolation method, one can apply the recursive extrapolation approach to obtain the solution on 
multi-time steps. }

\begin{figure}[H]
    \centering
    \includegraphics[width=0.32\linewidth]{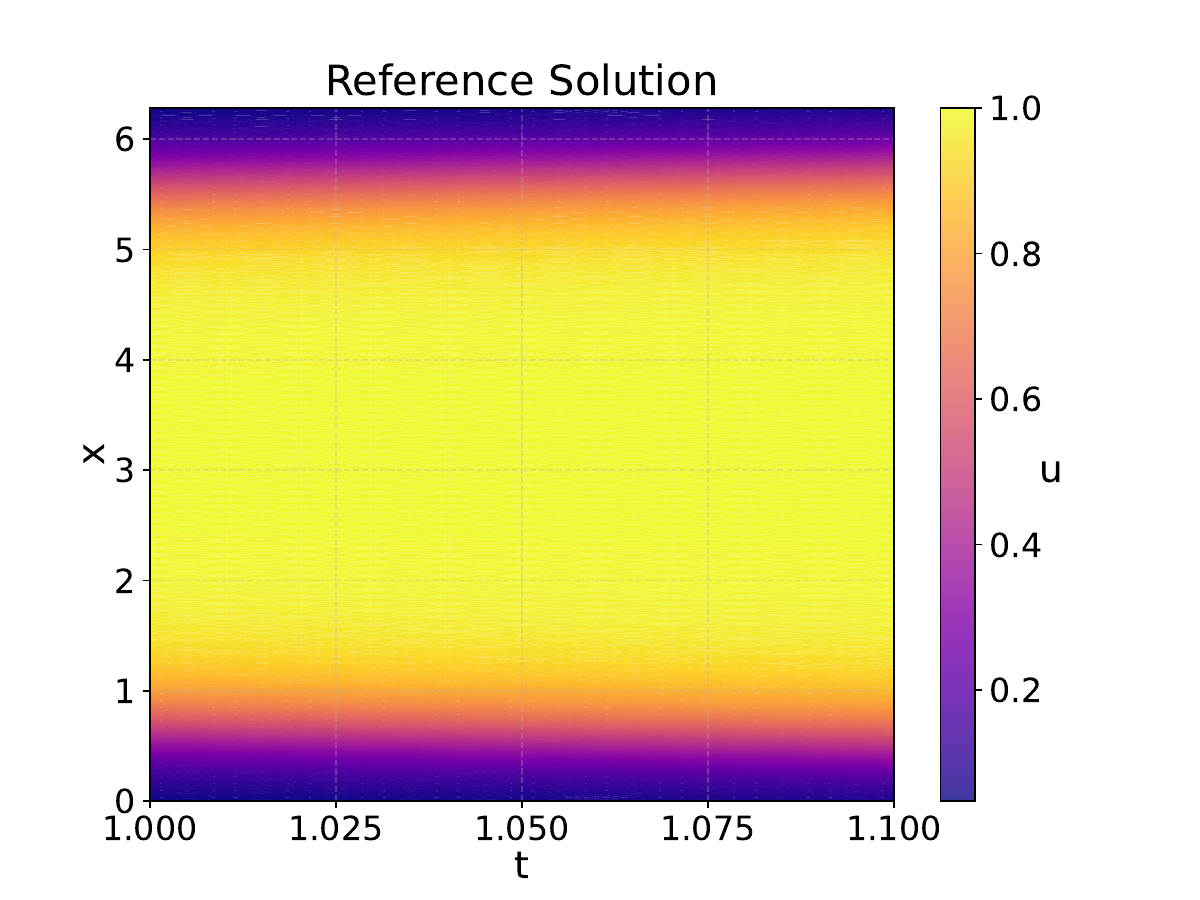}
    \includegraphics[width=0.32\linewidth]{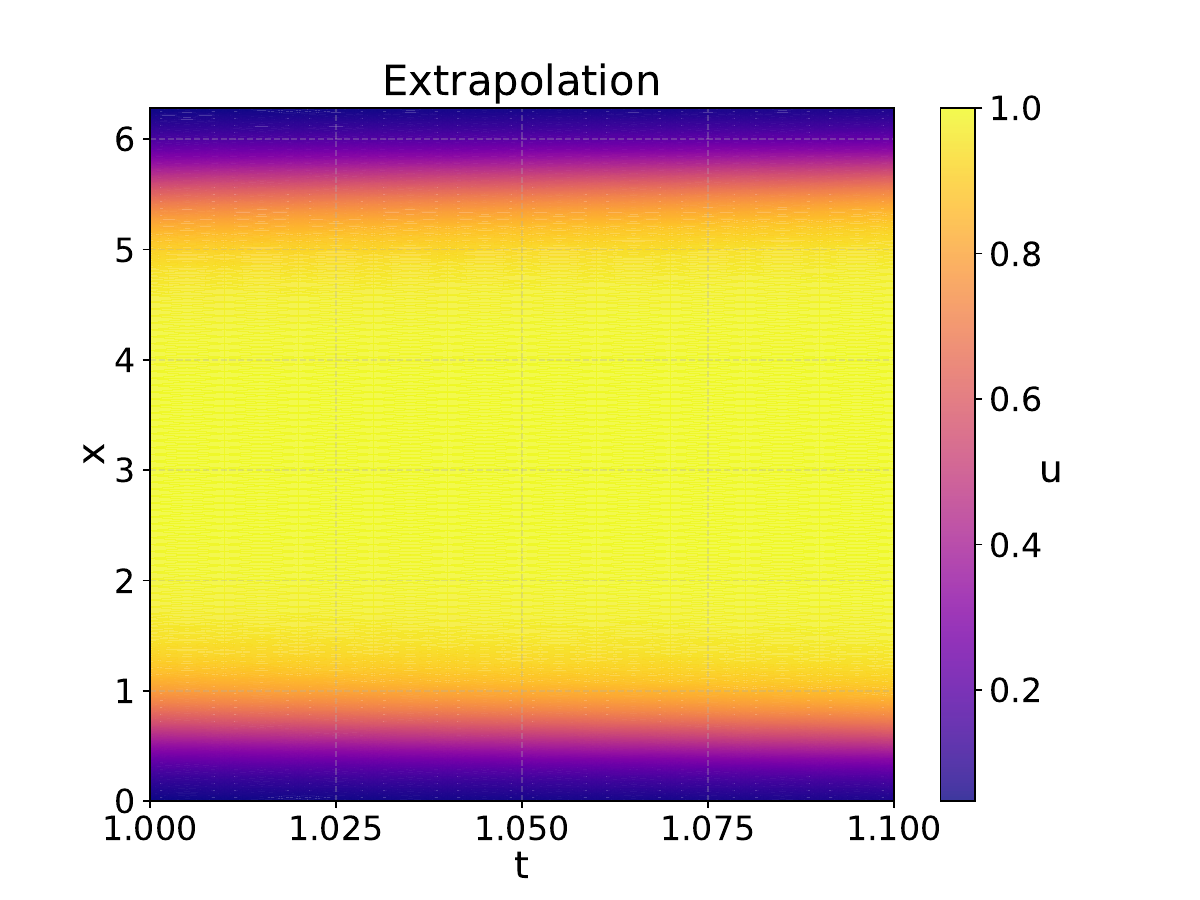}
    \includegraphics[width=0.32\linewidth]{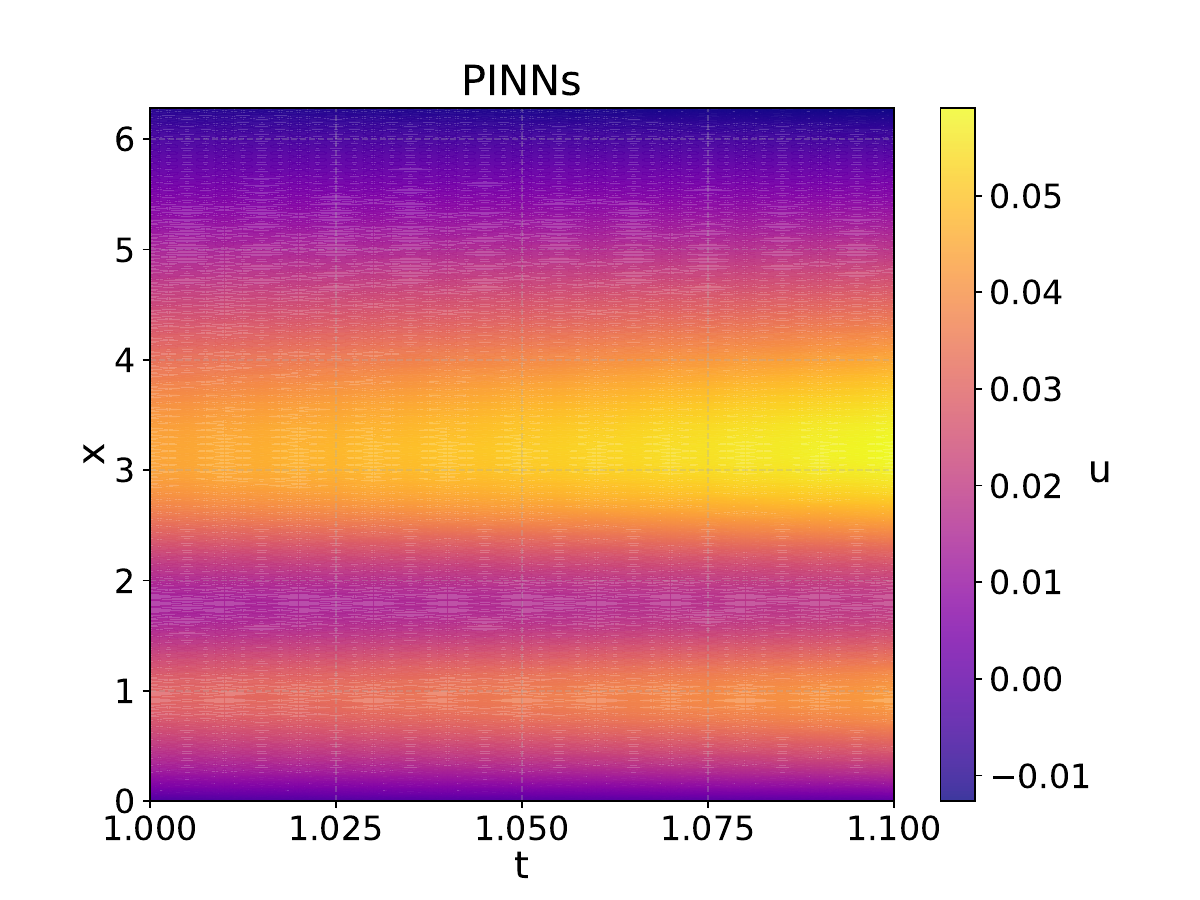}
    \includegraphics[width=0.32\linewidth]{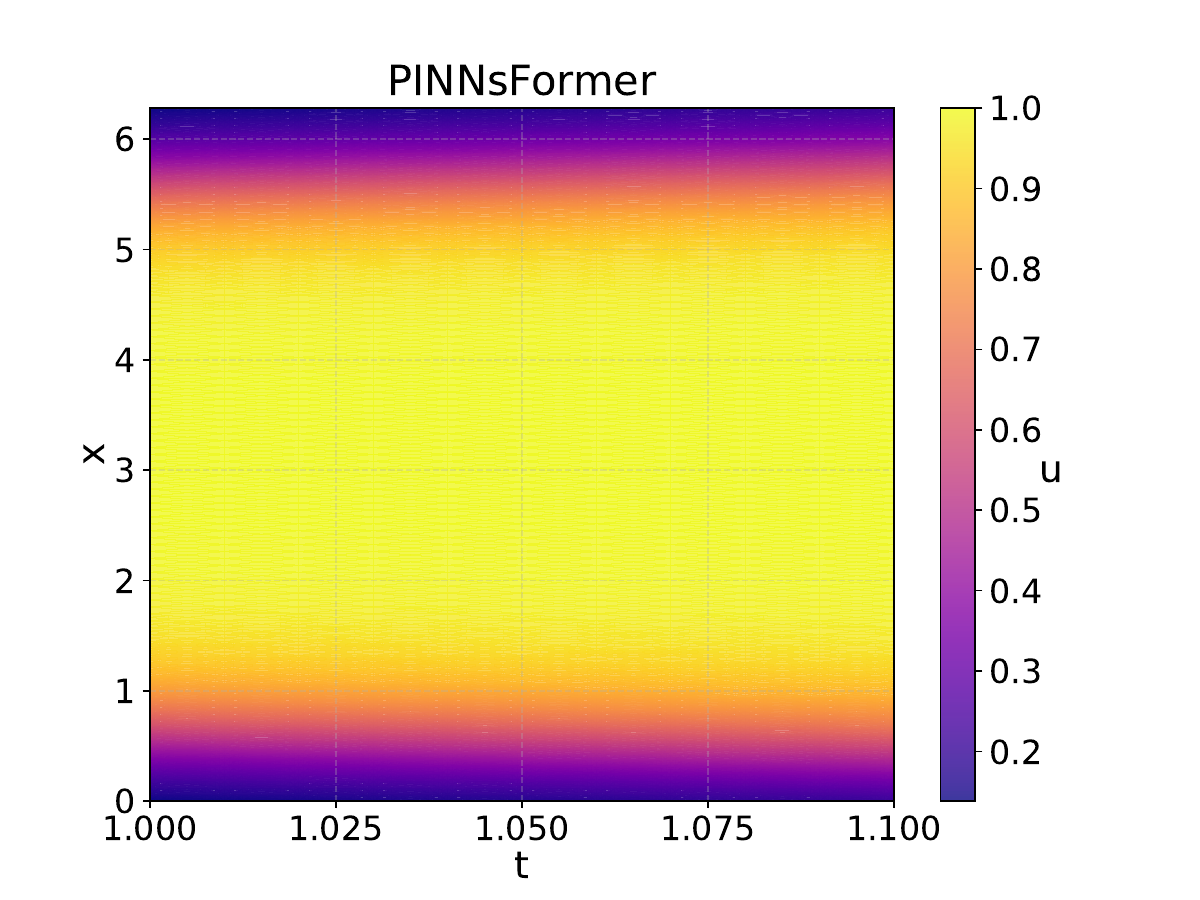}
    \includegraphics[width=0.32\linewidth]{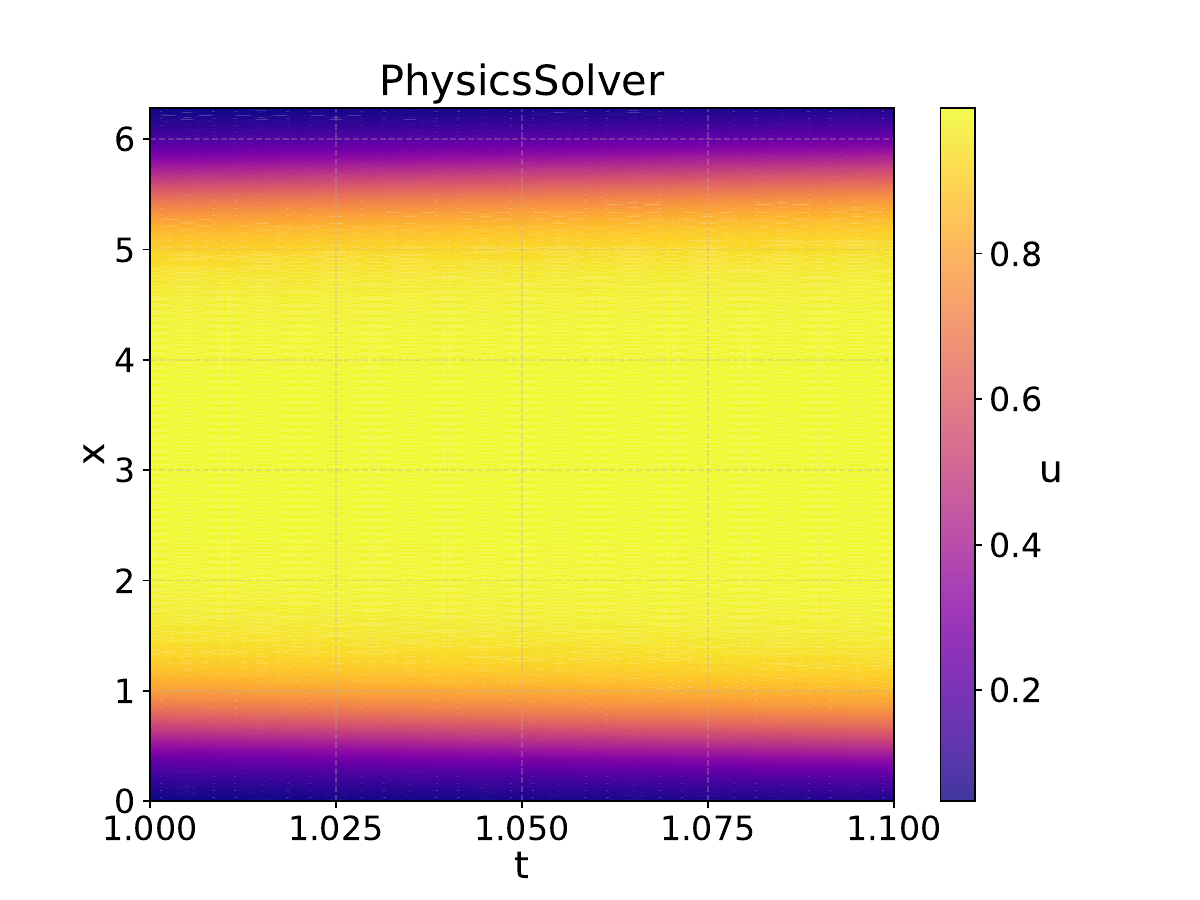}
    \caption{Forecasting problem with $\Delta t = 0.1$ across various methods for the reaction equation.}
    \label{fig:test5_t=10}
\end{figure}

\blue{First, we investigate the performance of different methods for $\Delta t=0.1$. The comparison results can be found in Figure \ref{fig:test5_t=10}. Table \ref{tab:test5_t=10} shows the detailed results which demonstrates that our proposed method outperforms other existing methods in terms of accuracy. We observe that our PhysicsSolver algorithm can achieve the relative $l^2$ error as small as $1.79\times 10^{-3}$. }

\begin{table}[H]
\vspace{-15pt}
    
	\vskip 0.1in
	\centering
	\begin{small}
		\begin{sc}
			\renewcommand{\multirowsetup}{\centering}
			\setlength{\tabcolsep}{5.5pt}
			\scalebox{1}{
			\begin{tabular}{l|c|c|c|c}
				\toprule
			    Method & Extrapolation & PINNs & PINNsFormer &PhysicsSolver \\
			    \midrule
                 Error &$5.015 \times 10^{-3}$& $9.673  \times 10^{-1}$  & $6.450 \times 10^{-2}$ & $\mathbf{1.793 \times 10^{-3}}$    \\
            \bottomrule
			\end{tabular}}
		\end{sc}
	\end{small}
    \caption{Relative $l^2$ errors for the forecasting problem with $\Delta t = 0.1$.}
	\label{tab:test5_t=10}
\end{table}

\begin{figure}[H]
    \centering
    \includegraphics[width=0.32\linewidth]{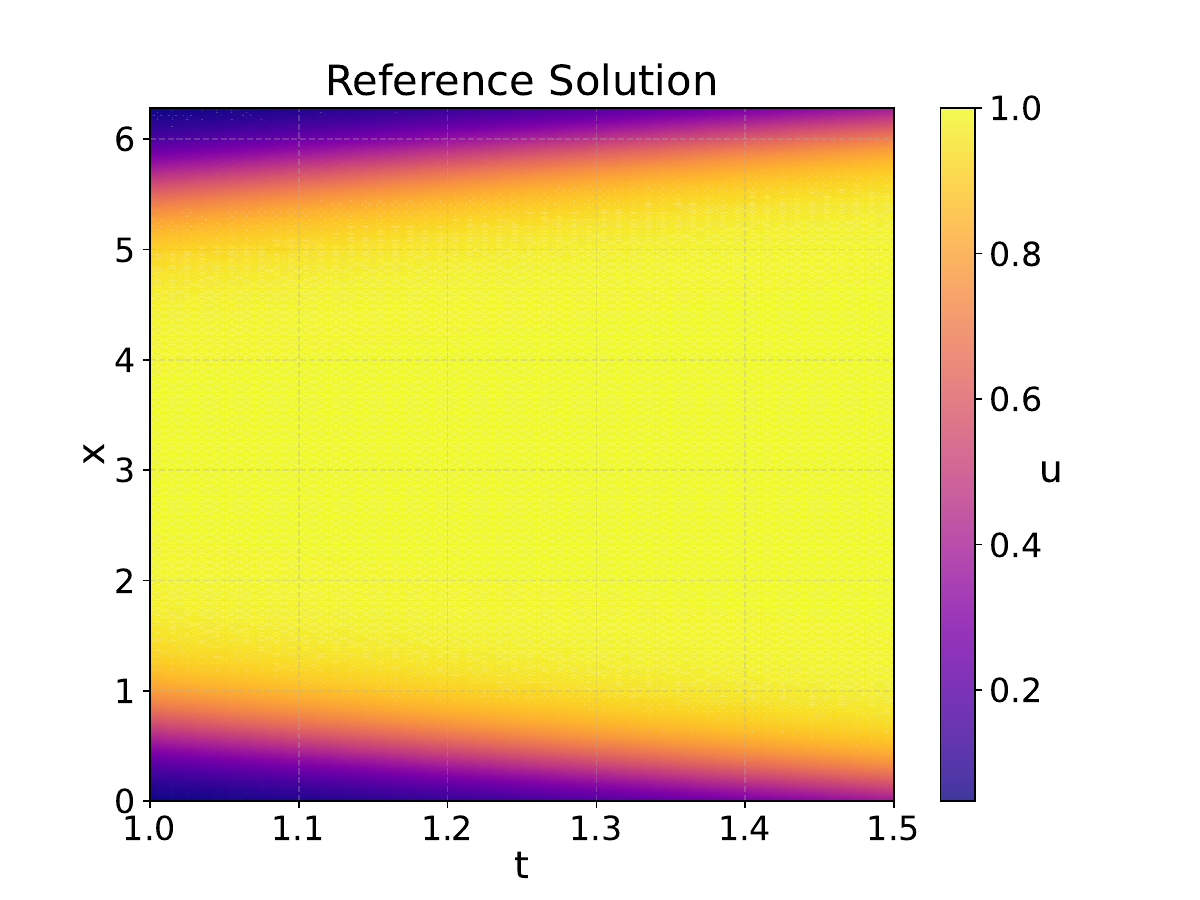}
    \includegraphics[width=0.32\linewidth]{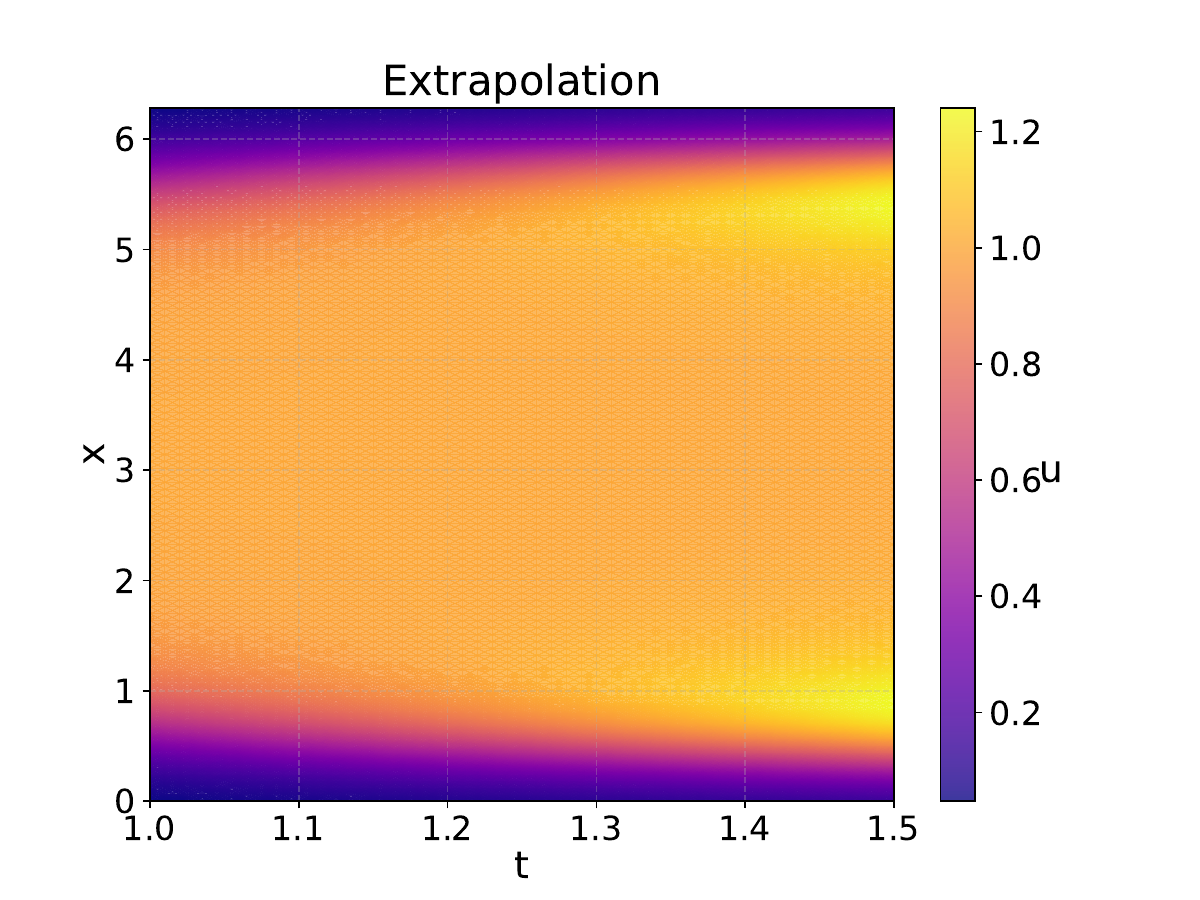}
    \includegraphics[width=0.32\linewidth]{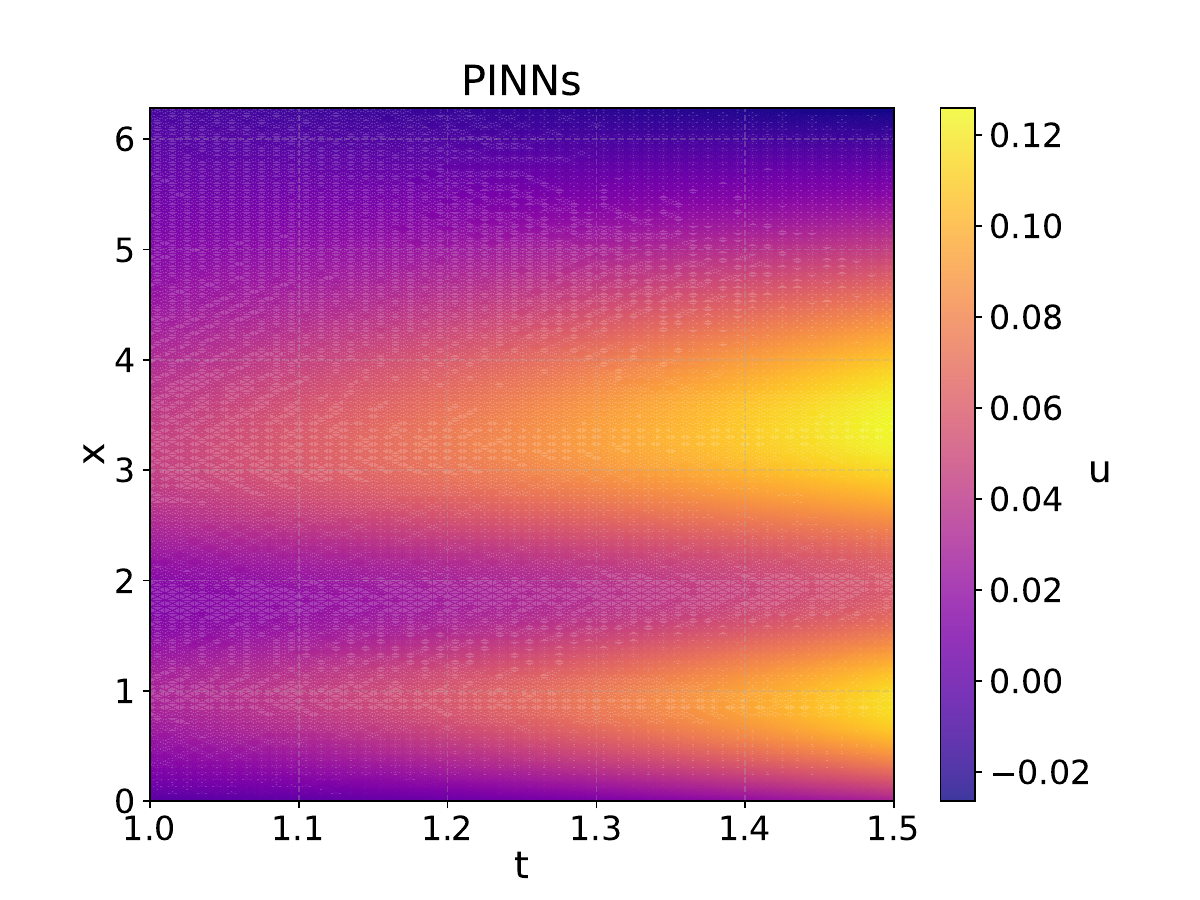}
    \includegraphics[width=0.32\linewidth]{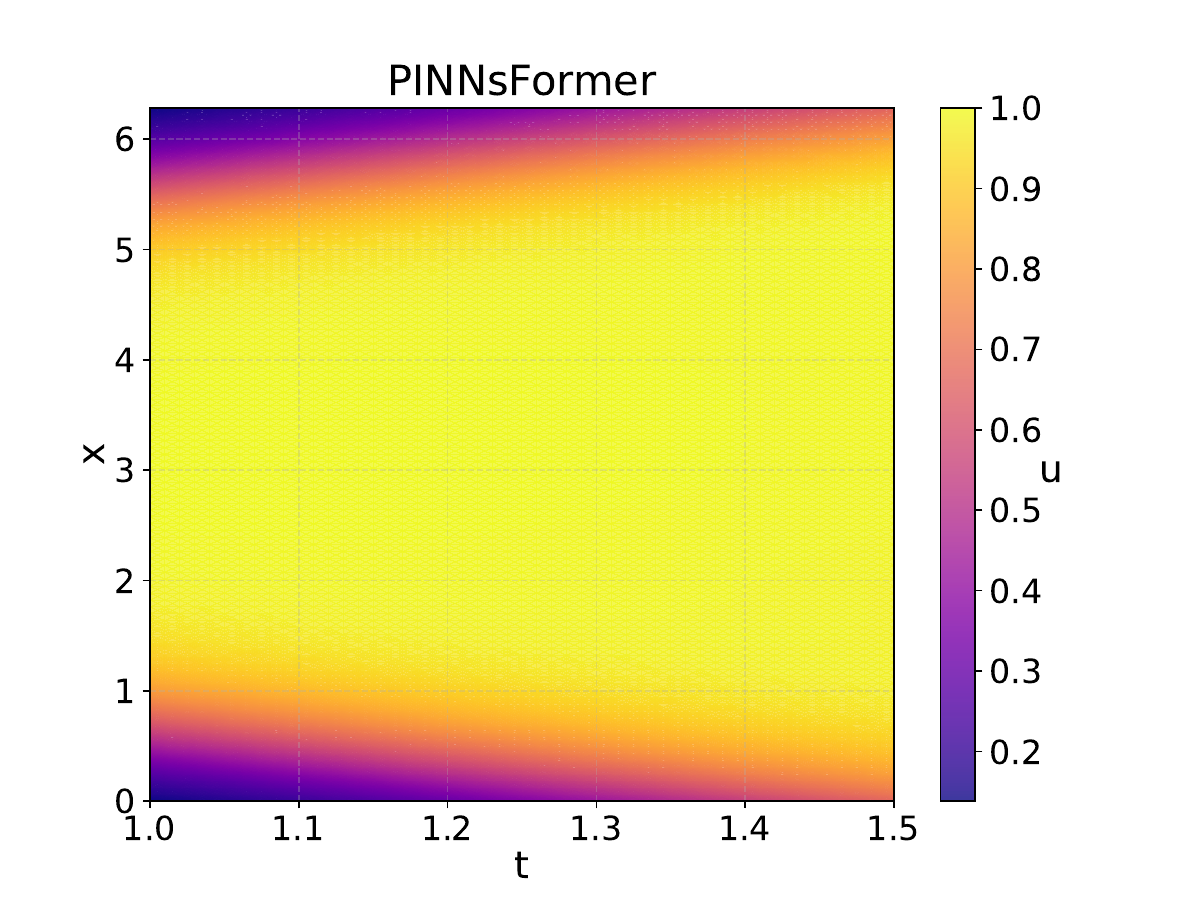}
    \includegraphics[width=0.32\linewidth]{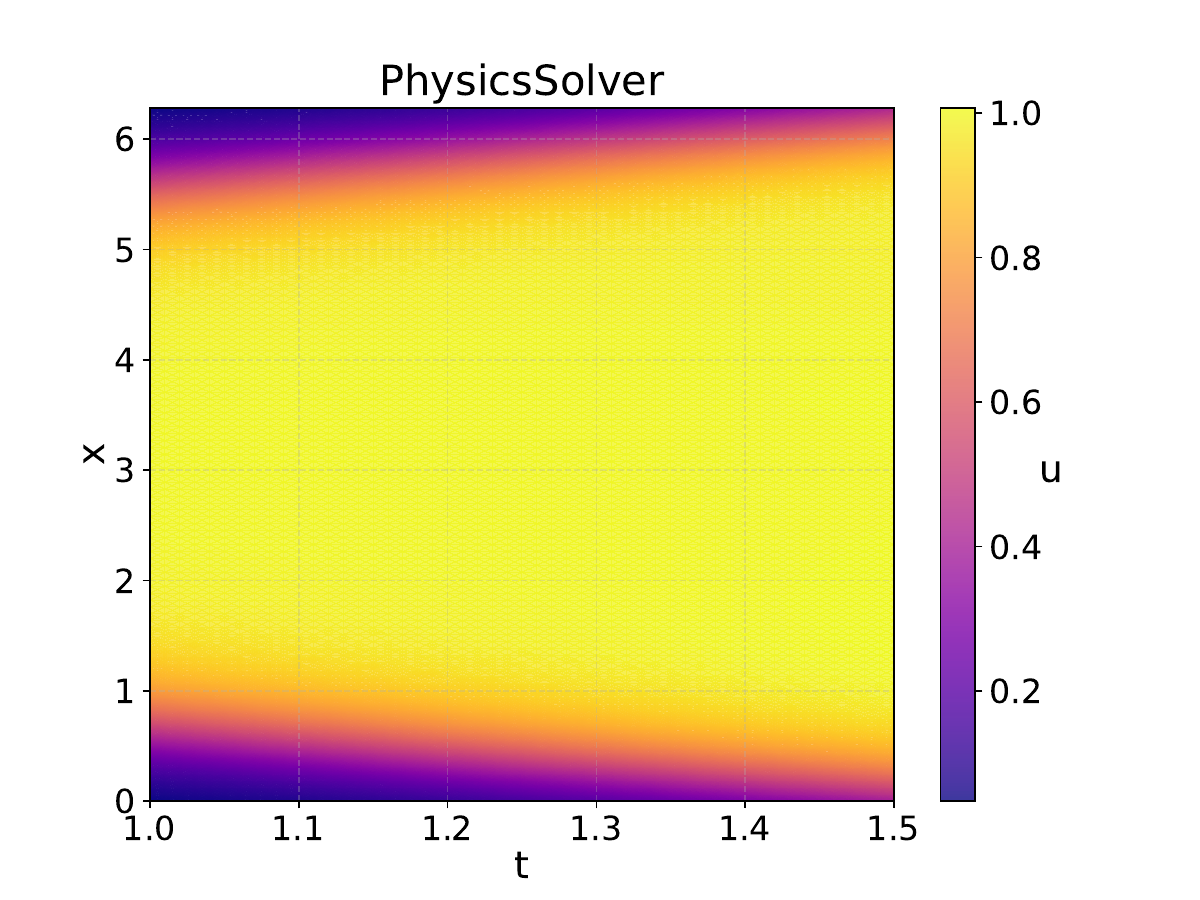}
    \caption{Forecasting problem with $\Delta t = 0.5$ across various methods for the reaction equation.}
    \label{fig:test5_t=50}
\end{figure}

\blue{Moreover, we investigate the performance of different methods for a larger time $\Delta t=0.5$. The comparison results are shown in Figure \ref{fig:test5_t=50} and Table \ref{tab:test5_t=50} for the relative $l^2$ errors of different approaches. We conclude that our PhysicsSolver significantly surpasses other methods in the long-term forecasting problem, reaching a level of accuracy as small as $\mathcal{O}(10^{-3})$. }

\begin{table}[H]
\vspace{-15pt}
    
	\vskip 0.1in
	\centering
	\begin{small}
		\begin{sc}
			\renewcommand{\multirowsetup}{\centering}
			\setlength{\tabcolsep}{5.5pt}
			\scalebox{1}{
			\begin{tabular}{l|c|c|c|c}
				\toprule
			    Method & Extrapolation & PINNs & PINNsFormer &PhysicsSolver \\
			    \midrule
                 Error &$7.639 \times 10^{-2}$& $9.462\times 10^{-1}$  & $7.785 \times 10^{-2}$ & $\mathbf{5.735 \times 10^{-3}}$    \\
            \bottomrule
			\end{tabular}}
		\end{sc}
	\end{small}
    \caption{Relative $l^2$ errors for the forecasting problem with $\Delta t = 0.5$.}
	\label{tab:test5_t=50}
\end{table}

\section{Conclusions and future work}
In this work, we propose an innovative transformer-based model called \textbf{PhysicsSolver}, which effectively solves both the forward and forecasting problems in PDE systems. Unlike previous approaches, our model can simultaneously learn intrinsic physical information and accurately predict future states by the ensemble of the physics system and the data system. We conduct extensive numerical experiments and the results show the superiority of our proposed method. 

\blue{\textit{Limitations and future work:} 
 First, we lack of theoretical analysis such as the convergence of predicted solution to reference solution in our work, as in most literature. Also, due to data limitations in certain examples, predicting PDEs with large time instants across more cases poses a significant challenge. In the future work, we plan to extend our method and study multi-step forecasting problem in broader applications. As time evolves, the error of extrapolation procedure may increase, requiring our current method to incorporate more constraints in the design of neural networks. Moreover, we intend to explore different variants of PhysicsSolver in tackling more complicated problems such as multi-scale kinetic models. }

\section*{Acknowledgement}
\label{sec:ack}
L.~Liu acknowledges the support by National Key R\&D Program of China (2021YFA1001200), Ministry of Science and Technology in China, Early Career Scheme (24301021) and General Research Fund (14303022 \& 14301423) funded by Research Grants Council of Hong Kong from 2021-2023.

\input{Appendix.tex}


\bibliographystyle{siam}
\bibliography{Submission_main_clean.bib}
\end{document}

%% file: Appendix.tex
\section*{Appendix}

\section*{Empirical Loss} 
The empirical loss for PINNs is:
\begin{equation}
\begin{aligned}
\mathcal{L}_\text{PINNs\_empirical} 
    =& w_r \sum_{i=1}^{N_\textit{r}} \|\mathcal{L}[\hat{u}(t_i, \boldsymbol{x}_i, \boldsymbol{v}_i)]-f(t_i, \boldsymbol{x}_i, \boldsymbol{v}_i)\|^2 \\
    &+
     w_b \sum_{i=1}^{N_b} \|\mathcal{B}[\hat{u}(t_i, \boldsymbol{x}_i, \boldsymbol{v}_i)]-g(t_i, \boldsymbol{x}_i, \boldsymbol{v}_i)\|^2 \\
    & + w_i\sum_{i=1}^{N_i} \|\mathcal{I}[\hat{u}(t_i, \boldsymbol{x}_i, \boldsymbol{v}_i)]-h(t, \boldsymbol{x}_i, \boldsymbol{v}_i)\|^2.
\end{aligned}
\label{empirical_loss_appendix}
\end{equation}
where $N_{r}$, $N_{b}$, and $N_{i}$ are the numbers of mesh grids used in computation for residual loss, boundary condition loss and initial condition loss, $w_r, w_b$, and $w_i$ are the corresponding weights, respectively.